\documentclass[review]{elsarticle}

\usepackage{lineno,hyperref}
\modulolinenumbers[5]

\journal{ }

\usepackage{commath}
\usepackage{dsfont}
\usepackage{lineno,hyperref}
\usepackage{amsmath,rotating}
\usepackage{subfigure}
\usepackage{psfrag}
\usepackage[para]{threeparttable} 
\usepackage[T1]{fontenc}
\usepackage{textcomp}
\usepackage{color}
\usepackage{placeins}
\usepackage{capt-of}
\usepackage{pstricks}
\usepackage{pst-node}
\usepackage{multido}
\usepackage{lineno}
\usepackage{hyperref}
\usepackage{tikz,pgfplots,tkz-fct,pgfplotstable}
\usepackage{siunitx}
\usepackage{framed,multirow}
\usepackage{fourier}
\tikzset{dashdot/.style={dash pattern=on .4pt off 3pt on 4pt off 3pt}}
\usetikzlibrary{decorations,decorations.text,backgrounds}
\usepackage{nicefrac}
\usetikzlibrary{%
    decorations.pathreplacing,%
    decorations.pathmorphing%
}
\usepackage{tkz-euclide}
\usetkzobj{all}
\usetikzlibrary{plotmarks}
\usetikzlibrary{shapes,arrows}
\usepgfplotslibrary{patchplots}

\usepackage{framed,multirow}

\usepackage{amssymb}
\usepackage{latexsym}

\usepackage{url}
\usepackage{xcolor}
\definecolor{newcolor}{rgb}{.8,.349,.1}

\newcommand\stent{\rule[-0.1ex]{0pt}{2.5ex}}

\newcommand\mrm[1]{\multirow{5}{*}[-0.5ex]{#1}}
\newcommand\mrmm[1]{\multirow{2}{*}[-0.5ex]{#1}}

\newcommand\lphi{$\;\; L_1(\bullet)\;\;\;$}
\newcommand\lphii{$\;\; L_2(\bullet)\;\;\;$}
\newcommand\lphiii{$\;\; L_{\infty}(\bullet)\;\;\;$}

\newcommand{\xGP}{-0.866, -0.339,  0.339, 0.866}
\newcommand{\xFV}{-0.75, -0.25,  0.25, 0.75}
\newcommand{\innerborders}{-0.5,0,0.5}

\begin{document}

\begin{frontmatter}

\title{A novel regularization strategy for the local discontinuous Galerkin method for level-set reinitialization}%
\author[1]{Fabian F\"oll\fnref{myfootnote}\corref{cor1}}
\ead{fabian.foell@iag.uni-stuttgart.de}
\author[1]{Christoph M\"uller\fnref{myfootnote}}
\ead{christoph.mueller@iag.uni-stuttgart.de}
\author[1]{Jonas Zeifang\fnref{myfootnote}}
\ead{jonas.zeifang@iag.uni-stuttgart.de}
\author[1]{Claus-Dieter Munz}
\ead{munz@iag.uni-stuttgart.de}
\address[1]{Institute of Aerodynamics and Gas Dynamics, University of Stuttgart, Pfaffenwaldring 21, 70569 Stuttgart, Germany}
\cortext[cor1]{Corresponding author}
\fntext[myfootnote]{These three authors contributed equally to the paper.}

\begin{abstract}
In this paper we propose a novel regularization strategy for the local discontinuous Galerkin method to solve the Hamilton-Jacobi equation in the context of level-set reinitialization. The novel regularization idea works in analogy to shock-capturing schemes for discontinuous Galerkin methods, which are based on finite volume sub-cells. In this spirit, the local discontinuous Galerkin method is combined with an upwind/downwind finite volume sub-cell discretization, which is applied in areas of low regularity. To ensure the applicability on unstructured meshes, the finite volume discretization is based on a least squares approach. 
\end{abstract}

\end{frontmatter}

\section{Introduction} \label{sec:Introduction}
	The level-set (LS) method, introduced by Osher and Sethian~\cite{osher198812}, is widely used in simulations of evolving interfaces. Instead of explicitly prescribing the shape of the interface, it is implicitly defined as the zero of the level-set function. Some recent developments and applications of the level-set method are summarized by Gibou et al.~\cite{gibou2018}. Although initially defined as a signed distance function with respect to the interface, the level-set field loses this property during integration in time due to velocity fields which occur, e.g. simulating multiphase flows~\cite{adalsteinsson1999,chang1996}. For the calculation of derivatives of the level-set field, such as the normal vector and the curvature, it is beneficial to retain the signed distance property. Several approaches exist for the reinitialization of the level-set field. They can be classified into explicitly reconstructing the interface, see e.g.~\cite{saye2014}, fast marching methods, as introduced in~\cite{sethian1996}, fast sweeping methods as proposed by~\cite{tsai2003} and flow based methods, see e.g.~\cite{sussman1994,karakus2016755}. In this work we focus on flow based methods, which are easy to parallelize~\cite{karakus2016755} and widely used in the application field of multiphase flows, see e.g.~\cite{fechter2015,winter2019,luo2019,gibou2019} for some recent applications. They rely on the solution of a Hamilton-Jacobi (HJ) type PDE
	\begin{align}
		\frac{\partial \phi}{\partial \tau} + H \left(\phi_x,\phi_y,\phi_z\right)= 0,
		\label{eq:HJ}
	\end{align} 
	with the Hamiltonian $H$. First order monotone finite difference (FD) schemes for solving HJ equations were developed by Crandall and Lions~\cite{crandall1984}. Later, Jiang and Peng~\cite{jiang2000} introduced a fifth order WENO scheme, which is limited to second order accuracy for solving the reinitialization equation. To overcome this, du Ch{\'e}n{\'e} et al.~\cite{chene2008} introduced a high order fix based on the work of Russo and Smereka~\cite{russo2000}. A main drawback of these methods is that an efficient parallelization is difficult due to their large stencil, especially on unstructured grids. Discontinuous Galerkin (DG) methods overcome this difficulty as they have a compact stencil although they are high order schemes. A DG scheme for solving HJ equations has been proposed by Hu and Shu~\cite{hu1999}. An extension of the first order monotone scheme in the DG sense has been introduced by Yan and Osher~\cite{yan2011232}. This local discontinuous Galerkin (LDG) scheme has arbitrary high order accuracy in space and reduces to the first order monotone scheme on Cartesian meshes for piecewise constant ansatz functions.\\
	Unfortunately, high order DG schemes tend to oscillate in the presence of discontinuities and kinks. They may occur during the reinitialization process or if the level-set function is cut off. To account for these instabilities, stabilization strategies for the LDG method have been proposed by Grooss and Hestaven~\cite{grooss20063406}, Mousavi~\cite{mousavi2014} and Karakus et al.~\cite{karakus2016755}. They applied an additional arbitrary numerical viscosity term to regularize the HJ equation and hence solved the equation
	\begin{align}
		\frac{\partial \phi}{\partial \tau} + H \left(\phi_x,\phi_y,\phi_z\right) - \nabla\cdot(\nu\nabla\phi) = 0.
		\label{eq:EikonalVisc}
	\end{align} 
  A drawback of the artificial viscosity regularization is the necessity to choose a suitable value for $\nu$. In~\cite{karakus2016755,grooss20063406,mousavi2014} arbitrary values, e.g. $\nu=0.0005,~0.001,~0.002$, were chosen. They turned out to have a sufficient stabilization effect.
  
  In this work, we present a new regularization strategy for the LDG approach of Yan and Osher~\cite{yan2011232}, which is based on the idea of finite volume sub-cells. This idea was introduced by Huerta et al.~\cite{huerta2012} and modified by Sonntag and Munz~\cite{sonntag2017,sonntag2017diss} as a shock-capturing method for a discontinuous Galerkin scheme, which is used to solve the Navier-Stokes equations. In contrast to their work, we use a least squares method to calculate the gradients for the finite volume sub-cells, since this allows the simulation of unstructured meshes.
  The main advantage of this new regularization strategy is, that it adds the numerical viscosity consistently, so that the solution converges to the viscosity solution of Eq.~\eqref{eq:HJ} and not of Eq.~\eqref{eq:EikonalVisc} (see \cite{sethian1996} for a definition of viscosity solutions). Moreover, no second order derivatives have to be calculated and the capability to fall back to a small stencil is introduced. The latter might be useful for the simulation of merging droplets or drop-wall interactions. Compared to finite volume and finite difference methods the benefits of the LDG scheme are preserved: arbitrary high order, handling of unstructured meshes and simple parallelization~\cite{cockDG2003,hesthaven2007}.
  
  The paper is structured as follows: In Sec.~\ref{sec:Governingequation} the relevant equations for the reinitialization are given. Then, the novel LDG based scheme is presented in Sec.~\ref{sec:LDG_LFV}. In Sec.~\ref{sec:Curvature}, a short note on the calculation of the gradients of the level-set field is made. Next, the capabilities of the novel method are illustrated with suitable test cases in Sec.~\ref{sec:Results} and a conclusion is drawn in Sec.~\ref{sec:Conclusion}.

\section{Governing equations}\label{sec:Governingequation}
\subsection{The level-set equation}\label{subsec:LSequation}
The level-set field $\phi$ is a signed distance function and implicitly describes the position of an interface located at the zero contour of the level-set field. Hence, it is defined as
\begin{align}
 	\phi(\boldsymbol{x})= \begin{cases}
 		- \mathcal{D}(\boldsymbol{x}) & \text{for}~\boldsymbol{x} \in \Omega^- \in \mathbb{R}^3, \\
 		+ \mathcal{D}(\boldsymbol{x}) & \text{for}~\boldsymbol{x} \in \Omega^+ \in \mathbb{R}^3, \\
 		0             & \text{for}~\boldsymbol{x} \in \Gamma,
 	\end{cases}
\end{align} 
where $\boldsymbol{x}=(x_1,x_2, x_3)^\text{T}=(x,y,z)^\text{T}$ and $\mathcal{D}(\boldsymbol{x})$ is the normal distance to the zero of the level-set field. The whole simulation domain is $\overline{\Omega} = \overline{\Omega^-}\cup \overline{\Omega^+}$ with $\Omega^- \cap \Omega^+ = \emptyset$ and the interface is $\Gamma = \overline{\Omega^-}\cap \overline{\Omega^+} $.
The normal vector on the interface $\boldsymbol{n}_{\Gamma}$ is related to the mean curvature $\kappa_{\Gamma}$ by
\begin{align}
	\kappa_{\Gamma} =  \nabla \cdot \boldsymbol{n}_{\Gamma}= \left. \nabla \cdot \frac{\nabla \phi}{| \nabla \phi |}\right \vert_{\Gamma} \quad \text{with}\quad \boldsymbol{n}_{\Gamma}= \left. \frac{\nabla \phi}{| \nabla \phi |}\right \vert_{\Gamma}.
	\label{eq:curvature}
\end{align}
The transport equation for the level-set field is
\begin{align}
	\frac{\partial \phi}{\partial t} + \boldsymbol{s} \cdot \nabla \left( \phi  \right) &=  0,
	\label{eq:levelset}
\end{align} 
with the level-set velocity field $\boldsymbol{s}$. 
The transport of the level-set field does not necessarily retain the signed distance property. This causes problems if derivatives of the level-set field have to be calculated. Therefore, a reinitialization procedure is applied to recover the signed distance property.

\subsection{The reinitialization equation} \label{subsec:HJequation}

Solving a HJ equation
\begin{align}
   \frac{\partial \phi}{\partial \tau} + H \left(\phi_x,\phi_y,\phi_z\right) &= 0,
   \label{eq:Eikonal}
\end{align} 
allows the reinitialization of the level-set function~\cite{sussman1994}, where $\tau$ is a pseudo time and $H \left(\phi_x,\phi_y,\phi_z\right)$ is the physical Hamiltonian defined by 
\begin{align}
  H\left(\phi_x,\phi_y,\phi_z\right) = \text{sign}(\phi)\left( | \nabla \phi| -1 \right).
   \label{eq:Hamilton}
\end{align} 
In the following we denote the vectors of level-set gradients in $x$, $y$, $z$ direction and in $\xi$, $\eta$, $\zeta$ direction as
\begin{align}
  \label{eq:lifting_vector_ref}
  \begin{split}
    \boldsymbol{\phi}^{\boldsymbol{x}}=\left[\phi_{x},\phi_{y},\phi_{z}\right]^\text{T}\quad \text{and} \quad  \boldsymbol{\phi}^{\boldsymbol{\xi}}=\left[\phi_{\xi},\phi_{\eta},\phi_{\zeta}\right]^\text{T}.
  \end{split}
\end{align}

\section{Regularization strategy for the local discontinuous Galerkin method} \label{sec:LDG_LFV}
\subsection{Semi-discrete form of the reinitialization procedure}\label{subsec:semi-discrete}
To solve the reinitialization equation (Eq.~\eqref{eq:Eikonal}) the physical Hamiltonian in Eq.~\eqref{eq:Hamilton} has to be approximated by a numerical Hamiltonian $H(\boldsymbol{\phi^{\boldsymbol{x}}})\approx\mathcal{H}(\boldsymbol{p},\boldsymbol{q}) $ using the vector of upwind gradients $\boldsymbol{p}\approx[\phi_{x},\phi_{y},\phi_{z} ]^\text{T}$ and the vector of downwind gradients $\boldsymbol{q}\approx[\phi_{x},\phi_{y},\phi_{z} ]^\text{T}$:
\begin{align}
\label{eq:lifting}
	\begin{split}
		p_{1}&-\phi_x^+=0,\quad  p_{2}-\phi_y^+=0,\quad  p_{3}-\phi_z^+=0, \\
		q_{1}&-\phi_x^-=0,\quad  q_{2}-\phi_y^-=0,\quad  q_{3}-\phi_z^-=0.
	\end{split}
\end{align} 
For the calculations in this work we use the Godunov Hamiltonian of the reinitialization equation
\begin{align}
	\label{eq:ldg_god}
	H(\boldsymbol{\phi^{\boldsymbol{x}}})&\approx\mathcal{H}^{God}(\boldsymbol{p},\boldsymbol{q}), \\
	&=\text{sign}(\phi)
	\begin{cases}
		\sqrt{\text{max}(a_1,b_1)+\text{max}(a_2,b_2)+\text{max}(a_3,b_3)-1}, & \text{if}\ \text{sign}(\phi) \leq 0, \\ \nonumber
		\sqrt{\text{max}(c_1,d_1)+\text{max}(c_2,d_2)+\text{max}(c_3,d_3)-1}, & \text{else},
	\end{cases} 
\end{align}
where 
\begin{align}
a_m=|p_m^+|^2,\quad b_m=|q_m^-|^2,\quad c_m=|p_m^-|^2,\quad d_m=|q_m^+|^2, \hspace{0.5cm} m=1,2,3,
	\label{eq:ldg_god_add1}
\end{align}
and
\begin{align}
f^+=\text{max}(f,0),\quad f^-=\text{min}(f,0).
	\label{eq:ldg_god_add2}
\end{align}
With this, the reinitialization procedure can be written in semi-discrete form as
\begin{align}
	\frac{\partial \phi}{\partial \tau }= -\mathcal{H}^{God}(\boldsymbol{p},\boldsymbol{q})\,.\label{eq:semi-discrete}
\end{align}

To further discretize Eq.~\eqref{eq:semi-discrete} Yan and Osher~\cite{yan2011232} introduced the LDG scheme, which is well suited in regions with smooth solutions but is unstable at discontinuities. In the following section we present a new regularization strategy for the LDG scheme. It is based on shock capturing strategies for non-smooth solutions of the Navier-Stokes equations, e.g. presented in~\cite{sonntag2017} and~\cite{sonntag2017diss}. The idea is to identify strong gradients and/or discontinuities in the solution with an indicator, see e.g.~\cite{persson2006}, and switch locally from the polynomial representation of the solution in the LDG scheme to a finite volume sub-cell representation. 
In the following, we first define the two building blocks of the solution: the LDG scheme and the corresponding finite volume sub-cell scheme. Next, the novel regularization strategy is presented, a numerical sign function is introduced and a short note on the time discretization is given.

\subsection{The lifting procedure and the local discontinuous Galerkin method } \label{subsec:LDG}
The LDG method is strongly related to the lifting procedure, see e.g.~\cite{bassi1997,arnold2002}. Hence, we first introduce the general lifting operator and in a second step the construction of the LDG scheme is described.
\subsubsection{The lifting procedure}
The lifting procedure has been introduced by Bassi and Rebay~\cite{bassi1997} for the approximation of the gradients that are required to calculate the parabolic fluxes of the Navier-Stokes equations.
The lifting operator in physical space in flux formulation is given by
\begin{align}
  \boldsymbol{d}  - \nabla_{\boldsymbol{x}} \cdot \underline{\boldsymbol{\Phi}} = 0,
 \label{eq:lifting_phys}
\end{align}
with the introduced lifting gradients $\boldsymbol{d}=\left[d^1,d^2,d^3\right]^\text{T}$ and the nabla operator in physical space $\nabla_{\boldsymbol{x}}\equiv\nabla=\left(\frac{\partial}{\partial x} ,\frac{\partial}{\partial y},\frac{\partial}{\partial z} \right)^\text{T}$.
Writing the level-set variable $\phi$ in a diagonal matrix simplifies the notation, $\underline{\boldsymbol{\Phi}}=\text{diag}\left([\phi,\phi,\phi]\right)$.
We subdivide the computational domain $\Omega$ with non-overlapping hexahedral elements $\Omega_e$ such that $\overline{\Omega} = \bigcup_e \overline{\Omega_e}$.
On arbitrary meshes, $\boldsymbol{x}(\boldsymbol{\xi})~\text{and}~\boldsymbol{\xi}(\boldsymbol{x})$ are the mappings between the coordinates of the physical space $\boldsymbol{x}$ and the reference space $\boldsymbol{\xi}=(\xi^1,\xi^2,\xi^3)^\text{T}=(\xi,\eta,\zeta)^\text{T}$.
Each element is mapped onto a reference cube element $\Omega^{ref}=[-1,1]^{3}$ using $\boldsymbol{\xi}(\boldsymbol{x})$.
The corresponding transformed equations to Eq.~\eqref{eq:lifting_phys} are given as
\begin{align}
  J \boldsymbol{d}  - \nabla_{\boldsymbol{\xi}} \cdot \underline{\boldsymbol{\Theta}} = 0\,,
 \label{eq:lifting_ref}
\end{align}
with the nabla operator in the reference space $\nabla_{\boldsymbol{\xi}}=\left(\frac{\partial}{\partial \xi} ,\frac{\partial}{\partial \eta},\frac{\partial}{\partial \zeta} \right)^\text{T}$ and the Jacobian
\begin{align}
  J=\boldsymbol{a}_i \cdot \left( \boldsymbol{a}_j \times \boldsymbol{a}_k \right)\,,  \hspace{0.5cm} (i,j,k)\in \mathbb{P}_3,
\end{align}
where $(i,j,k) \in \mathbb{P}_{3}:=\{(1,2,3),(2,3,1),(3,1,2)\}$ is defined as a \textit{cyclic permutation}. The covariant and contravariant basis vectors
\begin{align}
  \boldsymbol{a}_m:=\frac{\partial\boldsymbol{x}}{\partial \xi^m} \, \quad \text{and} \quad \hspace{0.5cm} \boldsymbol{a}^m:=\frac{\partial\boldsymbol{\xi}}{\partial x_m} \,, \hspace{0.5cm} m=1,2,3,
  \label{eq:metrics}
\end{align}
are defined for each element $\Omega^{ref}$ in the direction of the Cartesian coordinates $m$ in the reference plane. The corresponding transformed flux is
\begin{align}
  \label{eq:phi_ref}
 \underline{\boldsymbol{\Theta}}=\left[\boldsymbol{\Theta}^1,\boldsymbol{\Theta}^2,\boldsymbol{\Theta}^3\right]^\text{T}\,,\quad \boldsymbol{\Theta}^{m}=J \boldsymbol{a}^m \cdot \underline{\boldsymbol{\Phi}}\,, \hspace{0.5cm} m=1,2,3.
\end{align}
Multiplying Eq.~\eqref{eq:lifting_ref} by a test function $\varphi(\boldsymbol{\xi})$ and then integrating in space leads to the variational form of the lifting operator 
\begin{align}
\int_{\Omega^{ref}} \left( J \boldsymbol{d}  - \nabla_{\boldsymbol{\xi}} \cdot \underline{\boldsymbol{\Theta}}  \right) \varphi(\boldsymbol{\xi}) \text{d} \boldsymbol{\xi} =0\,.
\label{eq:ldg_variational}
\end{align}
We obtain the weak formulation 
\begin{align}
\underbrace{
  \int_{\Omega^{ref}} J \boldsymbol{d}  \varphi(\boldsymbol{\xi}) \text{d} \boldsymbol{\xi}}%
_{a} + 
\underbrace{
  \int_{\Omega^{ref}} \underline{\boldsymbol{\Theta}} \cdot \left( \nabla_{\boldsymbol{\xi}} \varphi(\boldsymbol{\xi}) \right)\text{d} \boldsymbol{\xi}}%
  _{b} - 
\underbrace{
\oint_{\partial \Omega^{ref}}   \underline{\boldsymbol{\Theta}} \cdot \boldsymbol{n}   \varphi(\boldsymbol{\xi}) \text{d} \boldsymbol{s}}%
_{c}
= \boldsymbol{0} \,,
\label{eq:ldg_weak}
\end{align}
by applying the divergence theorem. Note that $\boldsymbol{n}=[n_1,n_2,n_3]^\text{T}=[n_x,n_y,n_z]^\text{T}$ is the outward pointing unit normal vector of the element surface. There are three contributing parts: (a) a volume integral of variables $\boldsymbol{d}$, (b) a volume integral of variables $\underline{\boldsymbol{\Theta}}$ and (c) a surface integral of the fluxes.
In each element, the solution and the flux are approximated by polynomials,
\begin{align}
  \boldsymbol{d}_h = \sum_{i,j,k=0}^{\mathcal{N}} \boldsymbol{\hat{d}}_{ijk} \psi_{ijk}(\boldsymbol{\xi})
	\quad \text{and} \quad
  \boldsymbol{\Theta}_h^m = \sum_{i,j,k=0}^{\mathcal{N}} \boldsymbol{\hat{\Theta}}_{ijk}^m \psi_{ijk}(\boldsymbol{\xi}) \,, \hspace{0.5cm} m=1,2,3,
	\label{eq:dg_discret_ldg}
\end{align}
 where the test and basis functions are chosen identical according to Galerkin's idea.
Generally we collocate our solution at Legendre-Gauss nodes. The basis functions are products of the one-dimensional Lagrange polynomials of degree $\mathcal{N}$.
The integrals are evaluated by using Gauss quadrature rules, which are uniquely linked to the chosen collocation point set.
We approximate the flux at the element surface with a numerical flux function $\boldsymbol{\mathcal{G}}(\boldsymbol{\phi}^{int},\boldsymbol{\phi}^{ext}) \approx \underline{\boldsymbol{\hat{\Theta}}}_h \cdot \boldsymbol{n}$. We denote $(\bullet)^{int}$ as values from inside of the element and $(\bullet)^{ext}$ as values from the outside.

\subsubsection{The local discontinuous Galerkin method for the Hamilton-Jacobi equation}
Yan and Osher~\cite{yan2011232} applied the LDG method to obtain high order accurate approximations for the up- and downwind gradients $\boldsymbol{p}$, $\boldsymbol{q}$.
The lifting procedure given by Eq.~\eqref{eq:ldg_weak} and Eq.~\eqref{eq:dg_discret_ldg} is applied and the flux at the element surface is approximated with the numerical upwind flux $\boldsymbol{\mathcal{G}}^+(\boldsymbol{\phi}^{int},\boldsymbol{\phi}^{ext}) = \boldsymbol{\phi^+}$ or the numerical downwind flux $\boldsymbol{\mathcal{G}}^-(\boldsymbol{\phi}^{int},\boldsymbol{\phi}^{ext}) = \boldsymbol{\phi^-}$. The components of those fluxes are
\begin{align}
\phi^+_i=
                \begin{cases}
                  \phi^{ext}, & \text{if}\ n_i \geq 0, \\
                  \phi^{int}, & \text{else},
                \end{cases}  \hspace{0.5cm} i=1,2,3,
\label{eq:ldg_up}
\end{align}
and
\begin{align}
	\phi^-_i=
	\begin{cases}
		\phi^{int}, & \text{if}\ n_i \geq 0, \\
		\phi^{ext}, & \text{else},
	\end{cases}  \hspace{0.5cm} i=1,2,3.
	\label{eq:ldg_down}
\end{align}
Finally, the LDG method in semi-discrete form is defined as
\begin{align}
  \frac{\partial \phi}{\partial \tau }= -\mathcal{H}^{God}(\boldsymbol{p^{LDG}},\boldsymbol{q^{LDG}})\,.
\label{eq:DG_hamiltan}
\end{align}
Thereby $\boldsymbol{p^{LDG}}\equiv\boldsymbol{p}$ and $\boldsymbol{q^{LDG}}\equiv\boldsymbol{q}$ are the LDG gradients. The temporal discretization is discussed in Sec.~\ref{subsec:Timedisc}.

\subsection{Finite volume approximation of Hamilton-Jacobi equation} \label{subsec:LFV}
The second building block of our scheme is a first order finite volume (FV) discretization of the HJ equation. We define the gradients of the FV representation as
\begin{align}
\label{eq:FV_gradient_up}
\boldsymbol{p^{FV}}=\left[p^{FV}_1,p^{FV}_2,p^{FV}_3\right]^\text{T} = \left[\boldsymbol{D}^{x,+}\phi,\boldsymbol{D}^{y,+}\phi,\boldsymbol{D}^{z,+}\phi\right]^\text{T}, \\
\boldsymbol{q^{FV}}=\left[q^{FV}_1,q^{FV}_2,q^{FV}_3\right]^\text{T} = \left[\boldsymbol{D}^{x,-}\phi,\boldsymbol{D}^{y,-}\phi,\boldsymbol{D}^{z,-}\phi\right]^\text{T},
\label{eq:FV_gradient_down}
\end{align}
with the upwind and downwind difference operators in $\bullet$-directions $\boldsymbol{D}^{\bullet,+}$ and $\boldsymbol{D}^{\bullet,-}$. If we assume a Cartesian mesh, the upwind and downwind difference operator in $x$-direction simplifies to 
\begin{align}
	\begin{split}
  (\boldsymbol{D}^{x,+}\phi)_{i,j,k}&= \frac{\phi_{i,j,k}-\phi_{i-1,j,k}}{\Delta x}, \\
  (\boldsymbol{D}^{x,-}\phi)_{i,j,k}&= \frac{\phi_{i+1,j,k}-\phi_{i,j,k}}{\Delta x}.
  	\end{split}
\label{eq:FV_Operator_updown}
\end{align}
The gradients in $y$- and $z$-direction can be calculated analogously.
However, on unstructured hexahedral meshes all element-local directions $\xi, \eta, \zeta$ have to be considered. Therefore, we calculate the gradients with a least squares method, which is first introduced in general and then applied to the calculation of the upwind and downwind gradients.
 
The least squares method is based on the solution of the equation system 
\begin{align}
  (\phi_{ext_m} - \phi_{int}) =  (x_{ext_m} - x_{int})  \cdot \phi_x +  (y_{ext_m} - y_{int})  \cdot \phi_y +  (z_{ext_m} - z_{int})  \cdot \phi_z,
\end{align}
with $m=1,\dots,M$.
A compact notation is given by
\begin{align}
  \llbracket\phi_m \rrbracket  = \llbracket x_m \rrbracket \cdot \phi_x + \llbracket y_m \rrbracket \cdot \phi_y + \llbracket z_m \rrbracket \cdot \phi_z,\quad\quad\quad m=1, \dots , M,
\end{align}
with the jumps of the level-set variable $\llbracket \phi_m \rrbracket$ and the distances of the barycenters between the current element and the neighboring elements $\llbracket x_m \rrbracket,~\llbracket y_m \rrbracket,~\llbracket z_m \rrbracket$. The variable $M$ represents the number of jumps that contribute to the least squares gradients. With this we define the vector of jumps of the level-set variable as
\begin{align}
  \llbracket  \boldsymbol{\phi} \rrbracket  =
  \begin{pmatrix} 
  \llbracket\phi_{m=1} \rrbracket,
  \llbracket\phi_{m=2} \rrbracket,
  \cdots, 
  \llbracket\phi_{m=M} \rrbracket
\end{pmatrix}^{\text{T}}
\end{align}
and the matrix of distances of the barycenters as
\begin{align}
  \underline{\boldsymbol{M}} =
  \begin{pmatrix} 
  \llbracket x_{m=1}\rrbracket  & \llbracket y_{m=1}\rrbracket  & \llbracket z_{m=1}\rrbracket  \\
  \llbracket x_{m=2}\rrbracket  & \llbracket y_{m=2}\rrbracket  & \llbracket z_{m=2}\rrbracket  \\
  \cdots           & \cdots           & \cdots           \\
  \llbracket x_{m=M}\rrbracket  & \llbracket y_{m=M}\rrbracket  & \llbracket z_{m=M}\rrbracket  \\
  \end{pmatrix}.
\end{align}
As $\underline{\boldsymbol{M}}$ is not necessarily a square matrix the least squares method summarizes to 
\begin{align}
   \begin{split}
   \label{eq:leastsquares}
   \underline{\boldsymbol{M}} \boldsymbol{\phi}^{\boldsymbol{x}}  &  = \llbracket \boldsymbol{\phi} \rrbracket , \\
   \boldsymbol{\phi}^{\boldsymbol{x}}                 &              = \left(\underline{\boldsymbol{M}}^{\text{T}} \underline{\boldsymbol{M}}\right)^{-1} \underline{\boldsymbol{M}}^{\text{T}} \llbracket \boldsymbol{\phi} \rrbracket ,  \\
   \boldsymbol{\phi}^{\boldsymbol{x}}                 &              = \underline{\boldsymbol{W}} \llbracket \boldsymbol{\phi} \rrbracket,
   \end{split}
\end{align}
with the least squares matrix $\underline{\boldsymbol{W}}$ and $M\ge3$.

The FV reinitialization procedure requires upwind and downwind least squares operators, which are calculated in each direction. This reduces the number of possible jumps to $M=3$. We choose the correct barycenters and values of the neighboring cells in analogy to the LDG fluxes (Eq.~\eqref{eq:ldg_up} and Eq.~\eqref{eq:ldg_down}) for the upwind and downwind approximation
\begin{align}
&(\bullet)^{+}_{i}=
\begin{cases}
(\bullet)^{ext}, & \text{if}\ n_i \geq 0, \\
(\bullet)^{int}, & \text{else},
\end{cases}  \hspace{0.5cm} i=1,2,3\quad\text{and}\quad(\bullet)=\phi,x,y,z
\label{eq:fv_up}
\end{align}
and
\begin{align}
&(\bullet)^{-}_{i}=
\begin{cases}
(\bullet)^{int}, & \text{if}\ n_i \geq 0, \\
(\bullet)^{ext}, & \text{else},
\end{cases}  \hspace{0.5cm} i=1,2,3\quad\text{and}\quad(\bullet)=\phi,x,y,z.          
\label{eq:fv_down}
\end{align}
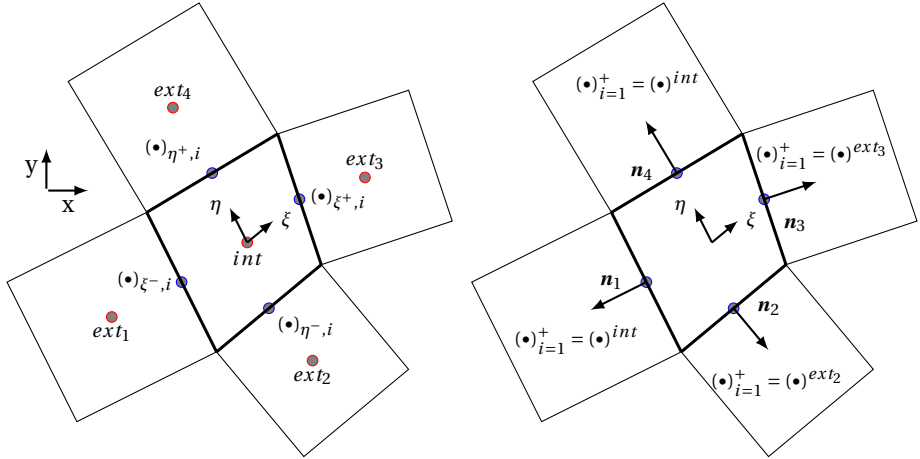
\begin{figure}[h!t]
	\begin{center}
		\begin{tikzpicture}[>=latex,baseline,scale=0.58]
		\tkzDefPoint(-0.4,0){A}
		\tkzDefPoint(2,2){B}
		\tkzDefPoint(1,5){C}
		\tkzDefPoint(-2,3.2){D}
		
		\tkzDefSquare(B,A)\tkzGetPoints{E}{F}
		\tkzDefSquare(C,B)\tkzGetPoints{G}{H}
		\tkzDefSquare(D,C)\tkzGetPoints{I}{J}
		\tkzDefSquare(A,D)\tkzGetPoints{K}{L}
		
		\tkzFillPolygon[draw,fill = white](A,B,F,E)
		\tkzFillPolygon[draw,fill = white](C,B,G,H)
		\tkzFillPolygon[draw,fill = white](C,D,J,I)
		\tkzFillPolygon[draw,fill = white](D,A,L,K)
		\tkzFillPolygon[draw,fill = white](A,B,C,D)
		
		\tkzInterLL(D,B)(A,C) 
		\tkzGetPoint{M}
		\tkzDrawPoint[size=10,color=red](M)
		
		\tkzInterLL(A,F)(B,E) 
		\tkzGetPoint{P}
		\tkzDrawPoint[size=10,color=red](P)
		
		\tkzInterLL(B,H)(C,G) 
		\tkzGetPoint{Q}
		\tkzDrawPoint[size=10,color=red](Q)
		
		\tkzInterLL(C,J)(D,I) 
		\tkzGetPoint{R}
		\tkzDrawPoint[size=10,color=red](R)
		
		\tkzInterLL(D,L)(A,K) 
		\tkzGetPoint{S}
		\tkzDrawPoint[size=10,color=red](S)
		
		\tkzInterLL(P,R)(Q,S) 
		\tkzGetPoint{T}
		
		\tkzDefMidPoint(A,B)
		\tkzGetPoint{AB}
		\tkzDrawPoint[size=10,color=blue](AB)
		
		\tkzDefMidPoint(B,C)
		\tkzGetPoint{BC}
		\tkzDrawPoint[size=10,color=blue](BC)
		
		\tkzDefMidPoint(C,D)
		\tkzGetPoint{CD}
		\tkzDrawPoint[size=10,color=blue](CD)
		
		\tkzDefMidPoint(D,A)
		\tkzGetPoint{DA}
		\tkzDrawPoint[size=10,color=blue](DA)
		
		\tkzInterLL(A,F)(B,E) 
		\tkzGetPoint{P}
		
		\tkzDefMidPoint(M,BC)
		\tkzGetPoint{MBC}
		
		\tkzDefMidPoint(M,CD)
		\tkzGetPoint{MCD}
		
		\draw [-latex, thick] (M) -- ($(MBC)$);
		\draw [-latex, thick] (M) -- ($(MCD)$);
		
		\tkzLabelPoint[right](MBC){\footnotesize{$\xi$}}
		\tkzLabelPoint[left](MCD){\footnotesize{$\eta$}}
		
		\tkzLabelPoint[below](M){\footnotesize{$int$}}
		\tkzLabelPoint[below](P){\footnotesize{$ext_2$}}
		\tkzLabelPoint[above](Q){\footnotesize{$ext_3$}}
		\tkzLabelPoint[above](R){\footnotesize{$ext_4$}}
		\tkzLabelPoint[below](S){\footnotesize{$ext_1$}}
		
		\tkzDrawSegment[-,line width=0.4mm](D,C)
		\tkzDrawSegment[-,line width=0.4mm](C,B)
		\tkzDrawSegment[-,line width=0.4mm](B,A)
		\tkzDrawSegment[-,line width=0.4mm](A,D)
		
		\tkzLabelPoint[below right](AB){\footnotesize{$(\bullet)_{\eta^-,i}$}}
		\tkzLabelPoint[right=0.03cm of BC](BC){\footnotesize{$(\bullet)_{\xi^+,i}$}}
		\tkzLabelPoint[above left=0.03cm of CD](CD){\footnotesize{$(\bullet)_{\eta^+,i}$}}
		\tkzLabelPoint[left](DA){\footnotesize{$(\bullet)_{\xi^-,i}$}}
		
		\node [inner sep=0pt] (cx) at (-4.3,4.7) {};
		\node [inner sep=0pt] (cy) at (-3.3,3.7) {};
		\node [inner sep=0pt] (o) at (-4.3,3.7) {};
		\draw [->,thick] (o) -- node[below] {x} (cy);
		\draw [->,thick] (o) -- node[left] {y} (cx);				
		\end{tikzpicture} \hspace{0.1cm}
		\begin{tikzpicture}[>=latex,baseline,scale=0.58]
		\tkzDefPoint(-0.4,0){A}
		\tkzDefPoint(2,2){B}
		\tkzDefPoint(1,5){C}
		\tkzDefPoint(-2,3.2){D}
		
		\tkzDefSquare(B,A)\tkzGetPoints{E}{F}
		\tkzDefSquare(C,B)\tkzGetPoints{G}{H}
		\tkzDefSquare(D,C)\tkzGetPoints{I}{J}
		\tkzDefSquare(A,D)\tkzGetPoints{K}{L}
		
		\tkzFillPolygon[draw,fill = white](A,B,F,E)
		\tkzFillPolygon[draw,fill = white](C,B,G,H)
		\tkzFillPolygon[draw,fill = white](C,D,J,I)
		\tkzFillPolygon[draw,fill = white](D,A,L,K)
		\tkzFillPolygon[draw,fill = white](A,B,C,D)
		
		\tkzInterLL(D,B)(A,C) 
		\tkzGetPoint{M}
		
		\tkzInterLL(A,F)(B,E) 
		\tkzGetPoint{P}
		
		\tkzInterLL(B,H)(C,G) 
		\tkzGetPoint{Q}
		
		\tkzInterLL(C,J)(D,I) 
		\tkzGetPoint{R}
		
		\tkzInterLL(D,L)(A,K) 
		\tkzGetPoint{S}
		
		\tkzInterLL(P,R)(Q,S) 
		\tkzGetPoint{T}
		
		\tkzDefMidPoint(A,B)
		\tkzGetPoint{AB}
		\tkzDrawPoint[size=10,color=blue](AB)
		
		\tkzDefMidPoint(B,C)
		\tkzGetPoint{BC}
		\tkzDrawPoint[size=10,color=blue](BC)
		
		\tkzDefMidPoint(C,D)
		\tkzGetPoint{CD}
		\tkzDrawPoint[size=10,color=blue](CD)
		
		\tkzDefMidPoint(D,A)
		\tkzGetPoint{DA}
		\tkzDrawPoint[size=10,color=blue](DA)
		
		\tkzInterLL(A,F)(B,E) 
		\tkzGetPoint{P}
		
		\tkzDefMidPoint(M,BC)
		\tkzGetPoint{MBC}
		
		\tkzDefMidPoint(M,CD)
		\tkzGetPoint{MCD}
		
		\draw [-latex, thick] (M) -- ($(MBC)$);
		\draw [-latex, thick] (M) -- ($(MCD)$);
		
		\tkzLabelPoint[right](MBC){\footnotesize{$\xi$}}
		\tkzLabelPoint[left](MCD){\footnotesize{$\eta$}}
		
		\draw [-latex, thick] (AB) -- ($(AB)!0.8!90:(A)$);
		\draw [-latex, thick] (BC) -- ($(BC)!0.8!90:(B)$);
		\draw [-latex, thick] (CD) -- ($(CD)!0.8!90:(C)$);
		\draw [-latex, thick] (DA) -- ($(DA)!0.8!90:(D)$);

		\tkzLabelPoint[right=0.2](AB){\footnotesize{$\boldsymbol{n}_2$}}
		\tkzLabelPoint[below right=0.2](BC){\footnotesize{$\boldsymbol{n}_3$}}
		\tkzLabelPoint[left=0.2](CD){\footnotesize{$\boldsymbol{n}_4$}}
		\tkzLabelPoint[left=0.2](DA){\footnotesize{$\boldsymbol{n}_1$}}
		
		\tkzLabelPoint[above](Q){\footnotesize{\,$(\bullet)^+_{i=1}=(\bullet)^{ext_3}\quad$}}
		\tkzLabelPoint[below](S){\footnotesize{$(\bullet)^+_{i=1}=(\bullet)^{int}$}}
		\tkzLabelPoint[below](P){\footnotesize{\quad$(\bullet)^+_{i=1}=(\bullet)^{ext_2}\quad$}}
		\tkzLabelPoint[above](R){\footnotesize{$(\bullet)^+_{i=1}=(\bullet)^{int}$}}
		
		\tkzDrawSegment[-,line width=0.4mm](D,C)
		\tkzDrawSegment[-,line width=0.4mm](C,B)
		\tkzDrawSegment[-,line width=0.4mm](B,A)
		\tkzDrawSegment[-,line width=0.4mm](A,D)
	
		\end{tikzpicture}
		
		\caption{Left: Nomenclature for the gradient approximation in two dimensions by the least squares approach; Right: Exemplary case for upwind gradient in x-direction (i=1).}
		\label{fig:leastsquare}
	\end{center}
\end{figure}
Details on the nomenclature can be found on the left of Fig.~\ref{fig:leastsquare}. An example for the upwind gradient in $x$-direction can be found on the right of Fig.~\ref{fig:leastsquare}.
With this, we define the vectors and matrices for the upwind and downwind least squares operators as
\begin{align}
  \boldsymbol{p^{FV}}   = \underline{\boldsymbol{W}}^+ \llbracket \boldsymbol{\phi^+} \rrbracket, \quad
  \boldsymbol{q^{FV}}   = \underline{\boldsymbol{W}}^- \llbracket \boldsymbol{\phi^-} \rrbracket,
  \label{eq:pq_FV}
\end{align}
with the jumps in each direction in the reference system
\begin{align}
  \llbracket      (\bullet)^\pm_{\xi,i}   \rrbracket =(\bullet)^\pm_{\xi_+,i}     -    (\bullet)^\pm_{\xi_-,i},    \quad
  \llbracket      (\bullet)^\pm_{\eta,i}  \rrbracket =(\bullet)^\pm_{\eta_+,i}    -    (\bullet)^\pm_{\eta_-,i},  \quad
  \llbracket      (\bullet)^\pm_{\zeta,i} \rrbracket =(\bullet)^\pm_{\zeta_+,i}   -    (\bullet)^\pm_{\zeta_-,i}.
\end{align}
Note that $\xi_\pm$, $\eta_\pm$ and $\zeta_\pm$ define values in positive and negative reference system directions. The matrices $\underline{\boldsymbol{W}}^+$ and $\underline{\boldsymbol{W}}^-$ in $x$, $y$ and $z$ direction are build-up with   
\begin{align}
  \begin{split}
    \underline{\boldsymbol{M}}^{x,\pm} & =
    \begin{pmatrix} 
    \llbracket x_{\xi  ,i=1}^{\pm}\rrbracket  & \llbracket y_{\xi  ,i=1}^{\pm}\rrbracket  & \llbracket z_{\xi  ,i=1}^{\pm}\rrbracket  \\[6pt]
    \llbracket x_{\eta ,i=1}^{\pm}\rrbracket  & \llbracket y_{\eta ,i=1}^{\pm}\rrbracket  & \llbracket z_{\eta ,i=1}^{\pm}\rrbracket  \\[6pt]
    \llbracket x_{\zeta,i=1}^{\pm}\rrbracket  & \llbracket y_{\zeta,i=1}^{\pm}\rrbracket  & \llbracket z_{\zeta,i=1}^{\pm}\rrbracket
    \end{pmatrix},\quad
    \underline{\boldsymbol{M}}^{y,\pm} =                                                                                           
    \begin{pmatrix}                                                                                                            
    \llbracket x_{\xi  ,i=2}^{\pm}\rrbracket  & \llbracket y_{\xi  ,i=2}^{\pm}\rrbracket  & \llbracket z_{\xi  ,i=2}^{\pm}\rrbracket  \\[6pt]
    \llbracket x_{\eta ,i=2}^{\pm}\rrbracket  & \llbracket y_{\eta ,i=2}^{\pm}\rrbracket  & \llbracket z_{\eta ,i=2}^{\pm}\rrbracket  \\[6pt]
    \llbracket x_{\zeta,i=2}^{\pm}\rrbracket  & \llbracket y_{\zeta,i=2}^{\pm}\rrbracket  & \llbracket z_{\zeta,i=2}^{\pm}\rrbracket
    \end{pmatrix},\\[6pt]
    \underline{\boldsymbol{M}}^{z,\pm} & =                          
    \begin{pmatrix}                                  
    \llbracket x_{\xi  ,i=3}^{\pm}\rrbracket  & \llbracket y_{\xi  ,i=3}^{\pm}\rrbracket  & \llbracket z_{\xi  ,i=3}^{\pm}\rrbracket  \\[6pt]
    \llbracket x_{\eta ,i=3}^{\pm}\rrbracket  & \llbracket y_{\eta ,i=3}^{\pm}\rrbracket  & \llbracket z_{\eta ,i=3}^{\pm}\rrbracket  \\[6pt]
    \llbracket x_{\zeta,i=3}^{\pm}\rrbracket  & \llbracket y_{\zeta,i=3}^{\pm}\rrbracket  & \llbracket z_{\zeta,i=3}^{\pm}\rrbracket
    \end{pmatrix}
  \end{split}
\end{align}
and result in 
\begin{align}                                                                                          
  \underline{\boldsymbol{W}}^\pm =                                                                                                    
  \begin{pmatrix}                                                                                                            
    \boldsymbol{D}^{x,\pm}, \boldsymbol{D}^{y,\pm}, \boldsymbol{D}^{z,\pm}  
\end{pmatrix}^{\text{T}}= 
  \begin{pmatrix}                                                                                                            
  W^{x,\pm}_{1,1}  & W^{x,\pm}_{1,2}  & W^{x,\pm}_{1,3}  \\[6pt]
  W^{y,\pm}_{2,1}  & W^{y,\pm}_{2,2}  & W^{y,\pm}_{2,3}  \\[6pt]
  W^{z,\pm}_{3,1}  & W^{z,\pm}_{3,2}  & W^{z,\pm}_{3,3}  
  \end{pmatrix}, 
\end{align}
according to Eqn. \eqref{eq:leastsquares}.
For the difference operator,~e.g. in x direction $\boldsymbol{D}^{x,\pm}$, only the corresponding contributions in $\underline{\boldsymbol{W}}^{x,\pm}$ are considered. For the gradients in y and z direction the procedure is applied analogously.
\\
In case of coordinate aligned elements $(\underline{\boldsymbol{M}}^{\bullet,\pm})^T\underline{\boldsymbol{M}}^{\bullet,\pm}$ can be singular and not invertible. This means that either too less or too much information about the direction-wise gradients is present to have a fully determined system. In such situations three different cases are distinguished according to the present zero columns and zero rows in $\underline{\boldsymbol{M}}^{\bullet,\pm}$:
\begin{itemize}
  \item After deleting the zero columns and rows $(\underline{\boldsymbol{M}}^{\bullet,\pm})^T\underline{\boldsymbol{M}}^{\bullet,\pm}$ is invertible: $\underline{\boldsymbol{W}}^{\bullet,\pm}$ is calculated according to Eq.~\eqref{eq:leastsquares} and filled up with the corresponding deleted zero columns and rows.
	\item After deleting the zero columns and rows only one information for each direction is present: The inverse can be directly calculated with Eq.~\eqref{eq:FV_Operator_updown} and filled up with the corresponding zero entries.
  \item After deleting the zero columns and rows $\underline{\boldsymbol{W}}^{\bullet,\pm}$ is underdetermined: $\underline{\boldsymbol{W}}^{\bullet,\pm}$ is calculated with the pseudo-inverse~\cite{bjoerck1996}
	\begin{align}
    \underline{\boldsymbol{W}}^{\bullet,\pm}=(\underline{\boldsymbol{M}}^{\bullet\pm})^T\left(\underline{\boldsymbol{M}}^{\bullet,\pm}(\underline{\boldsymbol{M}}^{\bullet\pm})^T\right)^{-1}.
	\end{align}
\end{itemize}
This procedure ensures that in each situation an optimal gradient (in the least squares sense) is calculated or is set to zero if the present data do not contain any information about the gradient in this direction. Note that, the evaluation of the matrices $\underline{\boldsymbol{M}}$ and $\underline{\boldsymbol{W}}$ is a preprocessing step. They do not need to be calculated during simulation run time. \\
Finally with Eq.~\eqref{eq:pq_FV}, we can define the low order semi-discrete scheme as
\begin{align}
	\frac{\partial \phi}{\partial \tau }= -\mathcal{H}^{God}(\boldsymbol{p^{FV}},\boldsymbol{q^{FV}})\,.
	\label{eq:FV_hamiltan}
\end{align}

\subsection{Novel regularization strategy for the local discontinuous Galerkin method} \label{subsec:Regularization}
The novel regularization strategy is based on the two previously defined building blocks:
\begin{enumerate}
  \item The LDG scheme, as described in Sec.~\ref{subsec:LDG}, is used in smooth regions with high regularity.
  \item The FV scheme, as described in Sec.~\ref{subsec:LFV}, is used in the presence of strong gradients or discontinuities with low regularity.
\end{enumerate}
Possible areas of low regularity in the level-set field are: (a) the edge of the narrow band where the level-set field is cut off, (b) areas with strong curvatures or (c) approaching level-set zeros, e.g. merging contours. The latter two cases may occur near the level-set zero. As a result, the combined method is both stable and high order accurate at the zero of the level-set function if possible.   

We introduce two different representations of the level-set field.
\begin{figure}[h!t]
\begin{center}

\begin{tikzpicture}[line cap=round,line join=round,x=1.6cm,y=1.6cm,axis/.style={->},>=latex]
 	  \foreach \a in {-1,1} {
 	      \draw (\a, -1) -- (\a, 1);
 	      \draw (-1, \a) -- (1, \a);
 	  }
 	  \foreach \x in \xGP {
 	      \foreach \y in \xGP {
 		  \draw [fill=black] (\x,\y) circle (0.03);
 	      }
 	  }
 	  \foreach \x in {-1,1} {
 	    \foreach \y in \xGP {
 		\draw[fill=white] (\x,\y) +(-0.03,-0.03) rectangle +(0.03,0.03);
 		\draw[fill=white] (\y,\x) +(-0.03,-0.03) rectangle +(0.03,0.03);
 	    }
 	  }
 	  \draw[axis] (0,0) node [anchor=south west] {$0$} -- (1.4,0) node [anchor=north east] {$\xi$};
 	  \draw[axis] (0,0) -- (0,1.4) node [anchor=north east] {$\eta$};
 	  \draw (1,0) node [anchor=south west] {$1$};
 	  \draw (0,1) node [anchor=south west] {$1$};
  \end{tikzpicture} \hspace{1cm}
\begin{tikzpicture}[line cap=round,line join=round,x=1.6cm,y=1.6cm,axis/.style={->},>=latex]
 	\foreach \a in {-1,1} {
 	    \draw (\a, -1) -- (\a, 1);
 	    \draw (-1, \a) -- (1, \a);
 	}
 	\foreach \x in \xFV {
 	    \foreach \y in \xFV {
 		\draw [fill=black] (\x,\y) circle (0.03);
 	    }
 	}
 	\foreach \x in {-1,1} {
 	  \foreach \y in \xFV {
 	      \draw[fill=white] (\x,\y) +(-0.03,-0.03) rectangle +(0.03,0.03);
 	      \draw[fill=white] (\y,\x) +(-0.03,-0.03) rectangle +(0.03,0.03);
 	  }
 	}
 	\foreach \y in {-2,0,2}{
 	  \draw[dashed] (-1.0,\y*0.25) -- (1.0,\y*0.25);
 	  \draw[dashed] (\y*0.25,-1) -- (\y*0.25,1);
 	}
 	\foreach \y in \innerborders {
 	  \foreach \x in \xFV {
 	      \draw[fill=white] (\x,\y) +(-0.03,-0.03) rectangle +(0.03,0.03);
 	      \draw[fill=white] (\y,\x) +(-0.03,-0.03) rectangle +(0.03,0.03);
 	  }
 	}
 	  \draw[axis] (0,0) node [anchor=south west] {$0$} -- (1.4,0) node [anchor=north east] {$\xi$};
 	  \draw[axis] (0,0) -- (0,1.4) node [anchor=north east] {$\eta$};
 	  \draw (1,0) node [anchor=south west] {$1$};
 	  \draw (0,1) node [anchor=south west] {$1$};
\end{tikzpicture}

   \caption{Idea of different solution representations of level-set field. Left: LDG representation on collocated nodes, right: FV sub-cell representation.}
   \label{fig:fvprojection}
\end{center}
\end{figure}
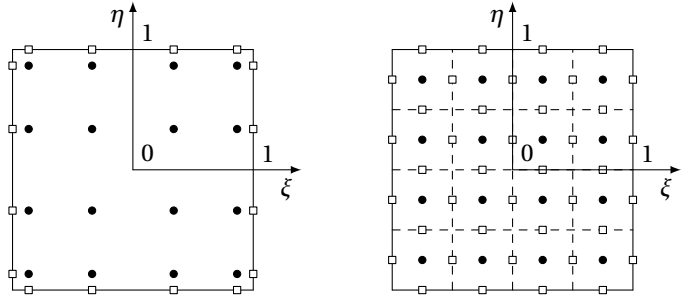
On the left side of Fig.~\ref{fig:fvprojection} the collocation points of the LDG representation are shown. The solution is represented by high order nodal polynomials. On the right side of Fig.~\ref{fig:fvprojection} the corresponding FV representation on an equidistant mesh is presented. The number of sub-cells coincides with the number of collocation points. During the reinitialization process the solution representation can be switched to the more robust FV approximation. We choose a conservative switch
\begin{align}
  \int_{\Omega^{ref}} J {\phi}  \text{d} \boldsymbol{\xi} \equiv \int_{\Omega^{ref}} J {\phi}_{LDG}\text{d} \boldsymbol{\xi} = \int_{\Omega^{ref}} J {\phi}_{FV}\text{d} \boldsymbol{\xi}
   \label{eq:FVDG}
\end{align} 
from the LDG to the FV representation and vice versa. This means a switch does not lead to a loss or growth of $\phi$ in general. The integrals are calculated exactly with a Gauss quadrature rule, which ensures conservation on the discrete level.
By applying the integration to each sub-cell a so called discrete projection matrix can be defined in one dimension, see Sonntag and Munz~\cite{sonntag2017}. The discrete projection matrix calculates the integral means of the FV sub-cells.
For the FV representation the scheme described in Sec.~\ref{subsec:LFV} is used, whereas for the polynomial representation the LDG scheme as described in Sec.~\ref{subsec:LDG} is equipped. A modal smoothness indicator in analogy to~\cite{huerta2012,persson2006} taken from~\cite{sonntag2017diss} is used to determine the suitable representation. The indicator is given by
\begin{align}
	\mathcal{S}=\text{log}_{10}\text{max}\Bigg\{ \frac{([\tilde{\phi}]_i^i,[\tilde{\phi}]_i^i)_{L_2}}{([\tilde{\phi}]_0^i,[\tilde{\phi}]_0^i)_{L_2}},i=\mathcal{N},\dots,\mathcal{N}-n\Bigg\},
\label{eq:Indicator}
\end{align}
with the $L_2$ inner product $(\cdot,\cdot)_{L_2}$. The truncation of modes between $a$ and $b$ is defined as
\begin{align}
	[\tilde{\phi}]_a^b=\sum_{i=l(a)-1}^{l(b)-1}\tilde{\phi}_i\mathcal{L}_i,
\end{align}
where $l(b)$ denotes the total number of the Legendre basis functions $\mathcal{L}_i$ for the polynomial degree $b$ and $n=1,2$. The transformation from nodal to modal space and vice versa can be achieved by the Vandermonde matrix $\underline{\boldsymbol{\mathcal{V}}}$ and its inverse $\underline{\boldsymbol{\mathcal{V}}}^{-1}$ by 
\begin{align}
  \boldsymbol{\tilde{\phi}}  = \underline{\boldsymbol{\mathcal{V}}}^{-1} {\boldsymbol{\hat{\phi}}} \quad \quad \text{and} \quad \quad \boldsymbol{\hat{\phi}} = \underline{\boldsymbol{\mathcal{V}}} {\boldsymbol{\tilde{\phi}}}\,.
\label{eq:modalnodalbasis}
\end{align}
Here $\boldsymbol{\tilde{\phi}}=[\tilde{\phi}_0,\dots,\tilde{\phi}_\mathcal{N}]^\text{T}$ is defined as the vector of modal values and $\boldsymbol{\hat{\phi}}=[\hat{\phi}_0,\dots,\hat{\phi}_\mathcal{N}]^\text{T}$ is defined as the vector of nodal values. Due to the tensor product structure, we can apply these operations line wise. So inherently, these are one-dimensional procedures. Note that we shift the zeroth mode of the level-set $\tilde{\phi}_0$ by one, in order to reduce the influence of the absolute value of the level-set on the indicator value. The indicator works in such a way that it evaluates the influence of the highest modes on the solution. Other indicators may be used, but it turned out that this indicator is suitable for the considered test cases.

The reinitialization procedure is done in the following way: First, the HJ-equation is solved by the LDG scheme (Eq.~\eqref{eq:DG_hamiltan}) on the whole domain. Second, the HJ-equation is solved by the FV scheme (Eq.~\eqref{eq:FV_hamiltan}) on the whole domain. Next, the indicator is evaluated on the whole domain (Eq.~\eqref{eq:Indicator}). A lower threshold $\mathcal{S}_{low}$ is defined to identify cells which remain in the LDG representation $\mathcal{S}\le\mathcal{S}_{low}$. 
In the same way an upper threshold $\mathcal{S}_{up}$ is defined to identify cells which require stabilization $\mathcal{S}\ge\mathcal{S}_{up}$.
For cells where $\mathcal{S}_{low}<\mathcal{S}<\mathcal{S}_{up}$ holds, the FV representation is switched to the nodal polynomial representation and the solutions are blended linearly. A considerable gain in efficiency can be obtained by evaluating the indicator in advance and only utilizing the LDG and the FV scheme in cells where they are required.

\subsection{Numerical sign function}
The sign function occurring in Eq.~\eqref{eq:Hamilton} and Eq.~\eqref{eq:ldg_god} is approximated by a smooth function
\begin{align}
\text{sign}(\phi)\approx\text{sgn}(\phi)=\phi/\sqrt{\phi^2 + \epsilon l_{ref}},
\end{align} 
with two characteristic length scales defined as
\begin{align}
l_{ref}=\min_{\forall \Omega_e \in \overline{\Omega}}(\sqrt[d]{V_e}) > 0 \quad \text{and} \quad \epsilon>0,
\end{align} 
where $V_e$ is the $d$-dimensional volume of the element $\Omega_e$. Its purpose is to smoothen the sign function over several cells. Here the parameter $l_{ref}$ explicitly considers the element size on arbitrary unstructured meshes, which may differ for each element. The parameter $\epsilon$ is then defined as a relative smoothing factor.

\subsection{Temporal discretization} \label{subsec:Timedisc}

A simple forward Euler time discretization of a HJ equation can be written as
\begin{align}
  \frac{\phi^{n+1}-\phi^{n}}{\Delta \tau }= -\mathcal{H}^{God}(\boldsymbol{p}^n,\boldsymbol{q}^n)\,,
\label{eq:Timedisc_HJ}
\end{align}
where the explicit time step restriction for the HJ equation  in one dimension is defined by the CFL condition 
\begin{align}
  \dfrac{\Delta t}{\Delta x} \lambda^{\text{1d}}  < \text{CFL} = 1 \quad \quad \text{and} \quad \quad \lambda^{\text{1d}} = \text{max}\abs{\text{sgn}(\phi)}_\Omega.
\end{align}
We approximate the time step in the multi-dimensional case on unstructured meshes by assuming the same eigenvalue as for the one dimensional case. Furthermore we redefine
\begin{align}
  \Delta x := \dfrac{2}{|\mathbf{a}^1|+|\mathbf{a}^2|+|\mathbf{a}^3|},
\end{align}
with $\mathbf{a}$ from Eq.~\eqref{eq:metrics} to account for stretched and deformed meshes. An extension of Eq.~\eqref{eq:Timedisc_HJ} to higher-order Runge-Kutta is straightforward. As the third order scheme from~\cite{Williamson1980} consists of consecutive applications of Eq.~\eqref{eq:Timedisc_HJ} we use this time integration scheme in some of the following applications.

\section{Gradient calculation for the level-set normals and the level-set curvature }\label{sec:Curvature}
For practical applications and for the evaluation of the accuracy of the scheme the gradient of the level-set field has to be calculated.
If the level-set has a polynomial representation, the simplest method for the calculation of the gradients is to directly derivate the polynomial basis of the level-set solution with
\begin{align}
  \nabla_{\boldsymbol{x}}\boldsymbol{\phi} \approx \sum_{i,j,k=0}^{\mathcal{N}}\hat{\boldsymbol{\phi}}_{ijk}\nabla_{\boldsymbol{x}}\boldsymbol{\psi}_{ijk}.
\end{align}
If the underlying solution is not smooth enough, jumps occur at element boundaries. In this case, the BR1 lifting procedure~\cite{bassi1997} is applied by using the numerical flux 
\begin{align}
  \underline{\boldsymbol{\hat{\Theta}}}_h \cdot \boldsymbol{n} \approx \boldsymbol{\mathcal{G}}^{\text{BR1}}(\boldsymbol{\phi}^{int},\boldsymbol{\phi}^{ext}) = \frac{1}{2}(\boldsymbol{\phi}^{ext}+\boldsymbol{\phi}^{int}),
\end{align}
in Eq.~\eqref{eq:ldg_weak}.
Other alternatives for the calculation of polynomial gradients are the $P_NP_M$ reconstruction method~\cite{dumbser2010} whose application to the level-set function was discussed in~\cite{fechter2015}, or other lifting procedures~\cite{fechter2015diss}.  

If the level-set has a FV representation, another method for the calculation of the gradients has to be considered. In this case we calculate gradients with the central least squares method Eq.~\eqref{eq:leastsquares}. Considering only the direct neighbors, $M=6$, the operator is build up with
\begin{align}
	\llbracket  \boldsymbol{\phi}^{\text{central}} \rrbracket  =
	\begin{pmatrix} 
		\llbracket\phi_{m=1} \rrbracket,
		\llbracket\phi_{m=2} \rrbracket,
		\llbracket\phi_{m=3} \rrbracket,
		\llbracket\phi_{m=4} \rrbracket,
		\llbracket\phi_{m=5} \rrbracket,
		\llbracket\phi_{m=6} \rrbracket
	\end{pmatrix}^{\text{T}},
\end{align}
and
\begin{align}
  \underline{\boldsymbol{M}}^{\text{central}} =
	\begin{pmatrix} 
		\llbracket x_{m=1}\rrbracket  & \llbracket x_{m=2}\rrbracket  & \llbracket x_{m=3}\rrbracket &\llbracket x_{m=4}\rrbracket  & \llbracket x_{m=5}\rrbracket  & \llbracket x_{m=6}\rrbracket  \\
		\llbracket y_{m=1}\rrbracket  & \llbracket y_{m=2}\rrbracket  & \llbracket y_{m=3}\rrbracket &\llbracket y_{m=4}\rrbracket  & \llbracket y_{m=5}\rrbracket  & \llbracket y_{m=6}\rrbracket  \\
		\llbracket z_{m=1}\rrbracket  & \llbracket z_{m=2}\rrbracket  & \llbracket z_{m=3}\rrbracket &\llbracket z_{m=4}\rrbracket  & \llbracket z_{m=5}\rrbracket  & \llbracket z_{m=6}\rrbracket  \\
	\end{pmatrix}^{\text{T}}.
\end{align}
Note that the central least squares matrices are always invertible. On structured, equidistant meshes this method reduces to the $2^{nd}$ order accurate central difference stencil. 
WENO methods~\cite{fechter2015diss} or finite difference approximations on structured grids~\cite{chene2008} are further options to calculate derivatives. \\
The described methods are used to calculate both the normal and the curvature of the level-set function according to Eq.~\eqref{eq:curvature}. For the calculation of the curvature the gradient operator is applied to the normalized field of the normal vector. The chosen variant to calculate both level-set normals and curvature is stated in each test case.

\section{Results} \label{sec:Results}
In this section the ability of the novel method to reinitialize disturbed level-set fields is shown. First, a convergence study for the two building blocks is considered. Second, the novel method is applied to benchmark test cases.
\subsection{Convergence tests} \label{subsec:Convergence}
%
In this section we provide results concerning  the convergence properties of the scheme. We use a test case in analogy to~\cite{chene2008} but choose a unit domain with $\Omega=[0,1]^2$. The level-set function is initially described by
\begin{align}
	\phi(x,y)=\text{exp}\left(10\sqrt{(x-0.5)^2+(y-0.5)^2}-2.313\right)-1.
\end{align}
The reinitialization procedure should give the signed distance function
\begin{align}
	\phi(x,y)=\sqrt{(x-0.5)^2+(y-0.5)^2}-0.2313,
\end{align}
with the corresponding curvature
\begin{align}
	\kappa=-\frac{1}{\sqrt{(x-0.5)^2+(y-0.5)^2}}.
\end{align}

\begin{figure}[t!]
	\begin{center}
		\begin{tikzpicture}[scale=0.99,>=latex]
			\node [inner sep=0pt] (center1) at (0,0)
				{\includegraphics[width=0.36\linewidth]{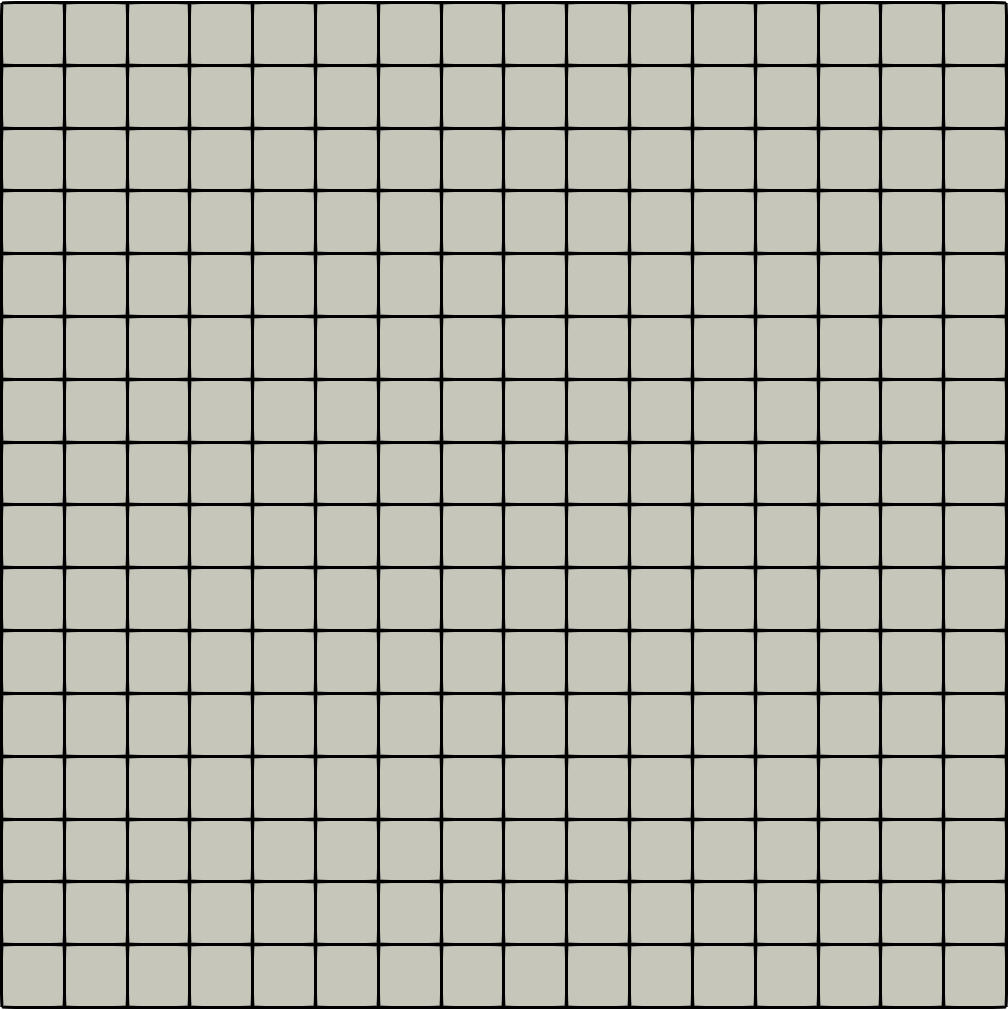}};
			\node [inner sep=0pt] (centerx1) at (-1.7cm,-2.7cm) {};
			\node [inner sep=0pt] (centery1) at (-2.7,-1.7cm) {};
			\node [inner sep=0pt] (origin1) at (-2.7cm,-2.7cm) {};
			\draw [->,thick] (origin1) -- node[left] {y} (centery1);
			\draw [->,thick] (origin1) -- node[below] {x} (centerx1);

			\node [inner sep=0pt] (center2) at (6,0)
				{\includegraphics[width=0.36\linewidth]{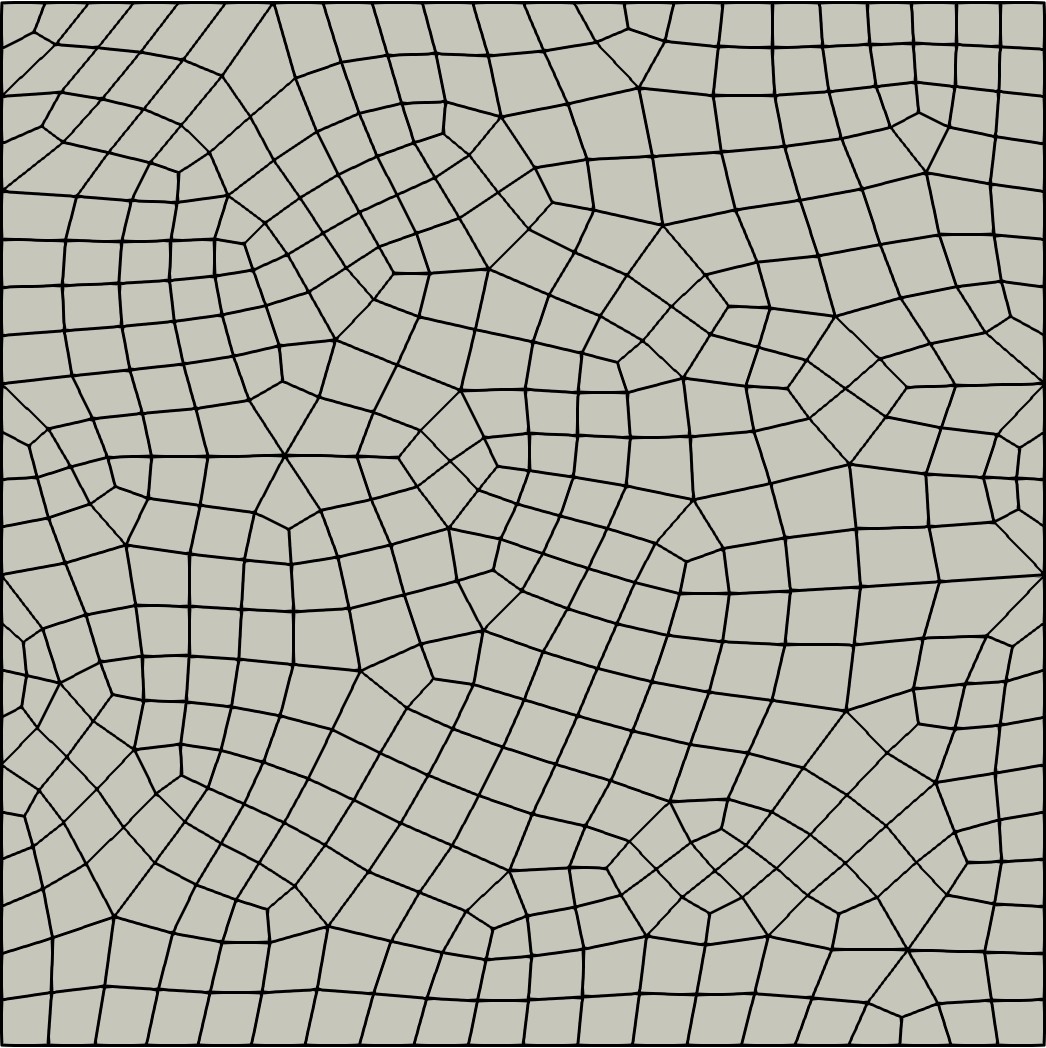}};
			\node [inner sep=0pt] (centerx2) at (4.3cm,-2.7cm) {};
			\node [inner sep=0pt] (centery2) at (3.3,-1.7cm) {};
			\node [inner sep=0pt] (origin2) at (3.3cm,-2.7cm) {};
			\draw [->,thick] (origin2) -- node[left] {y} (centery2);
			\draw [->,thick] (origin2) -- node[below] {x} (centerx2);				
		\end{tikzpicture}
		\caption{Exemplary meshes used for the convergence studies: structured mesh with 256 elements (left), unstructured mesh with 486 elements (right).}
		\label{fig:Meshes}
	\end{center}
\end{figure}
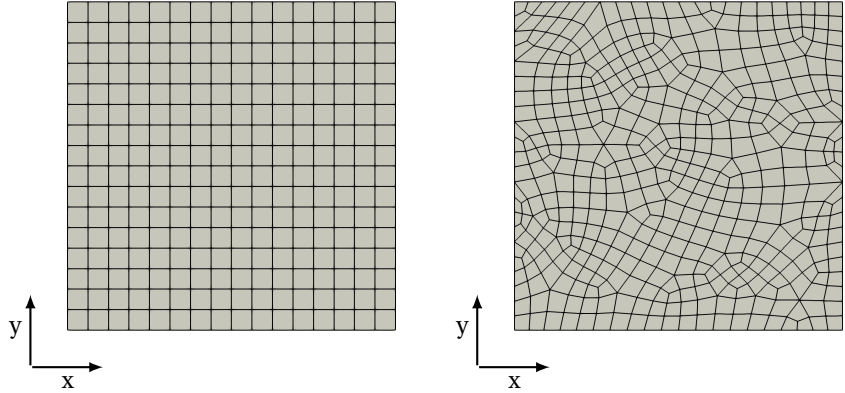

The simulations are performed with the forward Euler time discretization until a quasi-stationary solution is reached. In this context we define two termination criteria for the reinitialization process by
\begin{align}
||\phi^{n+1}-\phi^{n}||_{\infty} =  \Delta^n \leq \Xi \quad \vee \quad  \sum_{i=1}^N [\Delta^{n}- \Delta^{n-1} \le 0].
\label{eq:criteria}
\end{align}
Note that $\vee$ denotes a disjunction. The second criteria of Eq.~\eqref{eq:criteria} is a counter operator. For all calculations the parameters are chosen to $\Xi=1 \cdot 10^{-12}$ and $N=100$.
The convergence studies of the two building blocks are performed on Cartesian and unstructured meshes, as illustrated in Fig.~\ref{fig:Meshes}.
The $L_1,L_2,L_{\infty}$-errors are calculated without the kink in the center of the domain. Therefore, we exclude the innermost four cells near the kink. For the curvature $\kappa$ we also exclude the cells within $\{ 0.375<x_{\text{bary}}<0.625 \} \land \{ 0.375<y_{\text{bary}}<0.625 \}$ to avoid the influence of the kink. Note that the mathematical symbol $\land$ is a conjunction and the subscript $(\cdot)_{\text{bary}}$ refers to the barycenter of a cell. The curvature is calculated with the polynomial gradients and the central least squares gradients, respectively.
\begin{table}[b!]
\centering 
\begin{small}
  \begin{tabular}{|c|c|ll|ll|ll|}\hline
    \multirow{2}{*}{Variable } & mesh                              & \multicolumn{6}{|c|}{Polynomial degree of $\mathcal{N}=4$}  \\ \cline{3-8}
                               & \hspace{0.5em}level\hspace{0.5em} & \lphi                                                      & EOC  & \lphii   & EOC  & \lphiii  & EOC  \stent \\\hline\hline
    \mrm{level-set $\phi$}     & $h_0$                             & 4.53e-03                                                   & -    & 5.23e-03 & -    & 5.83e-03 & -    \stent \\
                               & $h_0/2$                           & 8.80e-05                                                   & 5.68 & 1.25e-04 & 5.37 & 1.61e-04 & 5.17 \stent \\
                               & $h_0/4$                           & 4.45e-06                                                   & 4.30 & 5.50e-06 & 4.51 & 8.75e-06 & 4.20 \stent \\
                               & $h_0/8$                           & 1.41e-07                                                   & 4.97 & 1.86e-07 & 4.88 & 4.41e-07 & 4.30 \stent \\\
                               & $h_0/16$                          & 5.41e-09                                                   & 4.71 & 1.02e-08 & 4.18 & 1.33e-07 & 1.72 \stent \\\cline{1-8}

    \mrm{curvature $\kappa$}   & $h_0$                             & 7.34e-01                                                   & -    & 9.71e-01 & -    & 1.20e+00 & -    \stent \\
                               & $h_0/2$                           & 5.97e-02                                                   & 3.62 & 7.75e-02 & 3.64 & 1.18e-01 & 3.34 \stent \\
                               & $h_0/4$                           & 7.20e-03                                                   & 3.05 & 1.02e-02 & 2.91 & 1.99e-02 & 2.56 \stent \\
                               & $h_0/8$                           & 1.06e-03                                                   & 2.76 & 1.61e-03 & 2.67 & 5.06e-03 & 1.97 \stent \\\
                               & $h_0/16$                          & 1.11e-04                                                   & 3.24 & 1.96e-04 & 3.03 & 7.69e-04 & 2.71 \stent \\\hline%
  \end{tabular}
\end{small}
\caption{h-convergence of the level-set and curvature field with $L_{\bullet}$ errors and $\mathcal{N}=4$ using the \text{\textbf{LDG scheme}} on \textbf{Cartesian meshes}, with $h_0 = 0.25$ and $\epsilon=50$.}
\label{tbl:LDG_Gauss_struct}
\end{table}

In Tbl.~\ref{tbl:LDG_Gauss_struct} the h-convergence of the level-set and curvature field, calculated with the LDG scheme on Cartesian grids, is given. For the level-set variable $\phi$ we observe EOC$\approx 5$, which is the theoretical order of convergence for a polynomial degree of $\mathcal{N}=4$. For the curvature variable $\kappa$ we observe EOC$\approx3$, which is the theoretical order of convergence. This is expected, since $\kappa$ is calculated from second derivatives of $\phi$. At the highest resolution, a flattening of the $L_{\infty}$-error in $\phi$ can be seen, which is also reported in~\cite{chene2008}. A possible explanation are inaccuracies due to roundoff errors.

In Tbl.~\ref{tbl:LDG_Gauss_unstruct} the h-convergence on unstructured meshes is shown. Similar to the results on the Cartesian meshes, both variables predominantly match the theoretical order of convergence. The inaccuracies are due to the fact that a refinement on unstructured grids cannot take place perfectly homogeneously.
The results of the h-convergence of the first order method are given in Tbl.~\ref{tbl:FV_conv} and Tbl.~\ref{tbl:FV_conv_unstructured}. Interestingly, for coarse resolutions, we observe EOC$\approx 2$ for the level-set variable, whereas for higher resolutions EOC$\approx 1$. This observation can be made on structured and unstructured meshes. The curvature shows no grid convergence as it can be expected for the low order method. The $L_\infty$ error of the curvature depends highly on the most deformed cell in the domain, which causes large differences of the $L_\infty$ error on different meshes.
\begin{table}[t!]
\centering 
\begin{small}
  \begin{tabular}{|c|c|ll|ll|ll|}\hline
   \multirow{2}{*}{Variable } & \#                                   & \multicolumn{6}{|c|}{Polynomial degree of $\mathcal{N}=4$}\\ \cline{3-8}
                              & \hspace{0.5em}elements\hspace{0.5em} & \lphi                                                      & EOC  & \lphii   & EOC  & \lphiii  & EOC   \stent \\\hline\hline
   \mrm{level-set $\phi$}     & $105$                                & 4.84e-05                                                   & -    & 6.94e-05 & -    & 1.75e-04 & -     \stent \\
                              & $291$                                & 4.74e-06                                                   & 4.56 & 6.92e-06 & 4.52 & 2.65e-05 & 3.70  \stent \\
                              & $718$                                & 5.60e-07                                                   & 4.72 & 7.51e-07 & 4.92 & 2.06e-06 & 5.66  \stent \\
                              & $1749$                               & 9.02e-08                                                   & 4.10 & 1.41e-07 & 3.74 & 5.86e-07 & 2.82  \stent \\\
                              & $4040$                               & 2.47e-08                                                   & 3.08 & 3.73e-08 & 3.17 & 1.36e-07 & 3.47  \stent \\\cline{1-8}
  
   \mrm{curvature $\kappa$}   & $105$                                & 1.28e-01                                                   & -    & 3.15e-01 & -    & 1.46e-00 & -     \stent \\
                              & $291$                                & 1.92e-02                                                   & 3.73 & 3.95e-02 & 4.07 & 1.79e-01 & 4.12  \stent \\
                              & $718$                                & 3.02e-03                                                   & 4.09 & 5.37e-03 & 4.41 & 2.46e-02 & 4.40  \stent \\
                              & $1749$                               & 1.29e-03                                                   & 1.89 & 2.35e-03 & 1.84 & 1.38e-02 & 1.29  \stent \\
                              & $4040$                               & 7.59e-04                                                   & 1.28 & 1.56e-03 & 0.97 & 1.17e-02 & 0.39 \stent \\\hline
  \end{tabular}
\end{small}
\caption{h-convergence of the level-set and curvature field with $L_{\bullet}$ errors and $\mathcal{N}=4$ using the \text{\textbf{LDG scheme}} on \textbf{unstructured meshes}, with $\epsilon=50$.}
\label{tbl:LDG_Gauss_unstruct}
\end{table}

In Tbl.~\ref{tbl:LDG_Gauss_pconv} a p-convergence study on structured grids is shown. It highlights a particular strength of the LDG method compared to WENO or finite difference methods: The LDG method is local and therefore spectral convergence can be achieved without a lot of implementation effort.

\begin{table}[t!]
	\centering 
	\begin{small}
		\begin{tabular}{|c|c|ll|ll|ll|}\hline
			 \multirow{2}{*}{Variable } & mesh  & \multicolumn{6}{|c|}{Effective polynomial degree of $\mathcal{N}=0$}\\ \cline{3-8}
			                          & \hspace{0.5em}level\hspace{0.5em} & \lphi    & EOC  & \lphii   & EOC  & \lphiii  & EOC  \stent \\\hline\hline
			 \mrm{level-set $\phi$}   & $h_0$                             & 2.30e-02 & -    & 2.44e-02 & -    & 4.02e-02 & -    \stent \\
			                          & $h_0/2$                           & 6.78e-03 & 1.76 & 7.66e-03 & 1.67 & 1.58e-02 & 1.35 \stent \\
			                          & $h_0/4$                           & 2.90e-03 & 1.22 & 3.40e-03 & 1.17 & 7.44e-03 & 1.08 \stent \\
			                          & $h_0/8$                           & 1.26e-03 & 1.20 & 1.53e-03 & 1.15 & 3.57e-03 & 1.06 \stent \\\
                                & $h_0/16$                          & 5.82e-04 & 1.12 & 7.20e-04 & 1.09 & 1.75e-03 & 1.03 \stent \\\cline{1-8}

			 \mrm{curvature $\kappa$} & $h_0$                             & 5.35e-01 & -    & 7.18e-01 & -    & 1.74e-00 & -      \stent \\
			                          & $h_0/2$                           & 3.82e-01 & 0.48 & 5.79e-01 & 0.31 & 2.40e-00 & -0.46  \stent \\
			                          & $h_0/4$                           & 1.52e-01 & 1.33 & 2.43e-01 & 1.25 & 2.21e-00 & 0.12   \stent \\
			                          & $h_0/8$                           & 1.13e-01 & 0.43 & 2.11e-01 & 0.21 & 2.53e-00 & -0.20  \stent \\\
			                          & $h_0/16$                          & 7.72e-02 & 0.55 & 1.96e-01 & 0.11 & 3.02e-00 & -0.25  \stent \\\hline%
		\end{tabular}
	\end{small}
  \caption{h-convergence of the level-set and curvature field with $L_{\bullet}$ errors using the \text{\textbf{FV method}} on \text{\textbf{Cartesian meshes}}, with $h_0 = 0.0625$ and $\epsilon=50$.}
	\label{tbl:FV_conv}
\end{table}
\begin{table}[b!]
	\centering 
	\begin{small}
		\begin{tabular}{|c|c|ll|ll|ll|}\hline
			 \multirow{2}{*}{Variable } & \#     & \multicolumn{6}{|c|}{Effective polynomial degree of $\mathcal{N}=0$}\\ \cline{3-8}
                                  & \hspace{0.5em}elements\hspace{0.5em} & \lphi    & EOC   & \lphii   & EOC   & \lphiii  & EOC  \stent \\\hline\hline
			                            & $744$                                & 1.74e-02 & -     & 1.93e-02 & -     & 3.46e-02 & -    \stent \\
       \mrmm{level-set $\phi$}    & $1812$                               & 7.97e-03 & 1.76  & 8.93e-03 & 1.73  & 2.07e-02 & 1.16 \stent \\
                                  & $3224$                               & 5.47e-03 & 1.31  & 6.05e-03 & 1.36  & 1.48e-02 & 1.16 \stent \\
                                  & $7500$                               & 3.15e-03 & 1.31  & 3.57e-03 & 1.24  & 8.28e-02 & 1.39 \stent \\\cline{1-8}

			                            & $744$                                & 5.79e-01 & -     & 7.98e-01 & -     & 4.08e-00 & -      \stent \\
       \mrmm{curvature $\kappa$}  & $1812$                               & 6.57e-01 & -0.28 & 9.79e-01 & -0.46 & 5.09e-00 & -0.50  \stent \\
                                  & $3224$                               & 7.51e-01 & -0.47 & 1.22e-00 & -0.76 & 1.43e+01 & -3.59  \stent \\
                                  & $7500$                               & 7.80e-01 & -0.09 & 1.26e-00 & -0.08 & 1.47e+01 & -0.07  \stent \\\hline%
		\end{tabular}
	\end{small}
  \caption{h-convergence of the level-set and curvature field with $L_{\bullet}$ errors using the \text{\textbf{FV method}} on \text{\textbf{unstructured meshes}}, with $\epsilon=50$.}
	\label{tbl:FV_conv_unstructured}
\end{table}
\begin{table}[t!]
\centering 
\begin{small}
  \begin{tabular}{|c|c|ll|ll|ll|}\hline
     \multirow{2}{*}{Variable } & polynomial    & \multicolumn{6}{|c|}{Number of elements $[n_x\times n_y] = [16 \times 16]$}\\ \cline{3-8}
                              & \hspace{0.5em}degree\hspace{0.5em} & \lphi    & EOC   & \lphii   & EOC   & \lphiii  & EOC          \stent\\\hline\hline
     \mrm{level-set $\phi$}   & $\mathcal{N}=0$                   & 2.59e-02 & -     & 2.72e-02 & -     & 3.89e-02 & -              \stent\\
                              & $\mathcal{N}=1$                   & 5.34e-03 & 2.27  & 5.81e-03 & 2.22  & 8.25e-03 & 2.23           \stent\\
                              & $\mathcal{N}=2$                             & 6.70e-04 & 5.12  & 8.55e-04 & 4.72  & 1.87e-03 & 3.65           \stent\\
                              & $\mathcal{N}=3$                   & 1.06e-04 & 6.40  & 1.29e-04 & 6.56  & 2.58e-04 & 6.88           \stent\\
                              & $\mathcal{N}=4$                   & 4.45e-06 & 14.20 & 5.50e-06 & 14.14 & 8.75e-06 & 15.17          \stent\\
                              & $\mathcal{N}=5$                   & 1.44e-07 & 18.79 & 1.92e-07 & 18.38 & 3.03e-07 & 18.43          \stent\\\cline{1-8}
     \mrm{curvature $\kappa$} & $\mathcal{N}=0$                   & 1.90e+00 & -     & 1.98e+00 & -     & 3.72e+00 & -     \stent\\
                              & $\mathcal{N}=1$                   & 1.83e+00 & 0.05  & 2.03e+00 & -0.03 & 3.61e+00 & 0.04   \stent\\
                              & $\mathcal{N}=2$                   & 6.60e-01 & 2.51  & 9.63e-01 & 1.83  & 2.44e+00 & 0.97   \stent\\
                              & $\mathcal{N}=3$                   & 1.33e-01 & 5.55  & 1.72e-01 & 5.98  & 3.29e-01 & 6.95   \stent\\
                              & $\mathcal{N}=4$                   & 7.20e-03 & 13.08 & 1.02e-02 & 12.63 & 1.99e-02 & 12.57  \stent\\
                              & $\mathcal{N}=5$                   & 5.02e-04 & 14.60 & 7.51e-04 & 14.35 & 1.47e-03 & 14.27  \stent\\\hline
  \end{tabular}
\end{small}
\caption{p-convergence of the level-set and curvature field with $L_{\bullet}$ errors using the \text{\textbf{LDG scheme}} on \textbf{Cartesian meshes}, with $\epsilon=50$. The case $\mathcal{N}=0$ was calculated with the FV scheme.}
\label{tbl:LDG_Gauss_pconv}
\end{table}

The results of the convergence studies proofed that the LDG method is high order accurate, even on unstructured meshes. Moreover, the results pointed out that a regularization should only be applied to elements away from the level-set zero position if possible, since a convergence of the curvature cannot be expected. It is important to note that if the underlying solution of the problem is a discontinuity, a high convergence order cannot be expected. Therefore it is permissible or even advisable to use a low order procedure even at the level-set zero position. Typical examples in multi-phase flow simulations are, e.g. merging droplet contours.
%
\subsection{Reinitialization of rectangular shaped level-set function} \label{subsec:Quad}
%
This test case is used to show the performance of the reinitialization method for discontinuous initial data. The discontinuity has a rectangular shape and is taken from~\cite{chene2008}. On the domain $\Omega=[-1,1]^2$, discretized with $33^2$ elements with a polynomial degree of $\mathcal{N}=4$, the level-set function is initially given by
\begin{align}
	\phi(x,y)=\begin{cases}
		+1, & \text{if}\ |x|\geq 0.5 ~ \text{or} ~ |y|\geq 0.5, \\
		-1, & \text{else}.
	\end{cases}
\end{align}
For the second test the initial data is rotated by $45$ degree in avoidance of aligning the discontinuity with the grid cells. In an application the level-set function is typically cut off by
\begin{align}
	\phi = \begin{cases}
		\phi, & \text{if}\ |\phi|\le \phi_\text{cut-off}, \\
		\text{sign}(\phi)\phi_\text{cut-off}, & \text{else},
	\end{cases}
\end{align}
to create a narrow band around the level-set zero. The reinitialization is only active in this band, which reduces the required amount of iterations. Here, we introduce a cut-off of the level-set function with $\phi_\text{cut-off}=0.25$. We observed favorable effects on the stability. The remaining parameters are: $\epsilon=20$, $CFL=0.5$, upper threshold of indicator $\mathcal{S}_{up}=-6.5$ and lower threshold of indicator is $\mathcal{S}_{low}=-7.5$ with $n=2$. The gradient of the level-set is calculated with the BR1 lifting procedure and the central least squares method.

\begin{figure}[b!]
    \centering
    \includegraphics[width=0.35\textwidth,trim={3cm 3cm 3cm 3cm},clip]{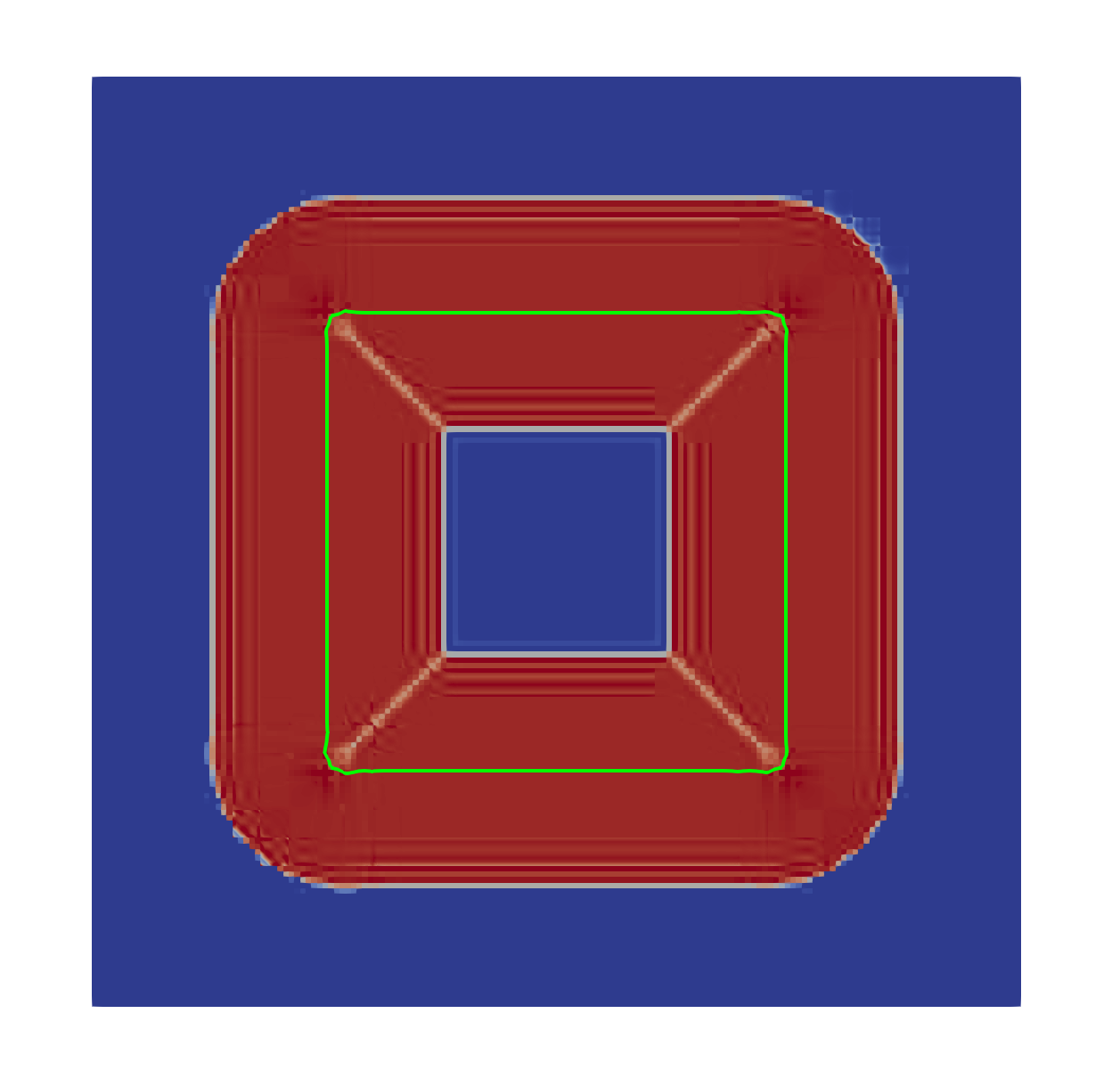}\hspace{0.5cm}
		\includegraphics[width=0.35\textwidth,trim={3cm 3cm 3cm 3cm},clip]{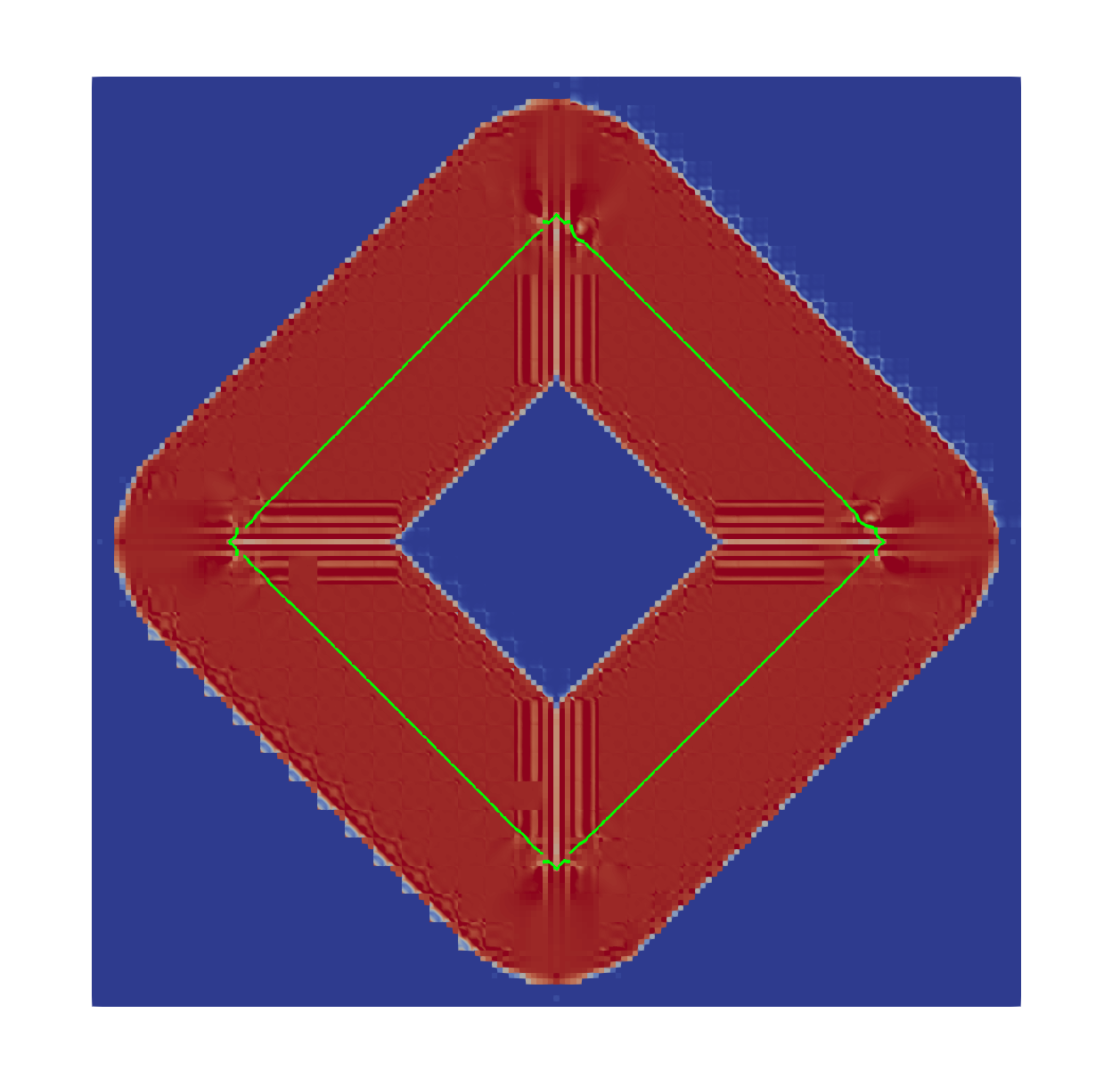}
    \vspace{0.2cm}
    \begin{tikzpicture}[line cap=round,line join=round,>=triangle 45,x=1.15cm,y=1.15cm,axis/.style={->}]
          \node at (0 ,-0.4) {\includegraphics[scale=0.10]{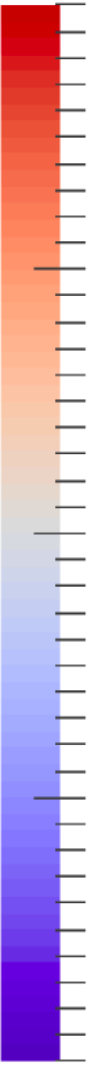}};
          \node at (0.0,1.6 ) {$~~~|\nabla \phi|~~~$};
          \node at (0.5,1.2 ) {$1.04$};
          \node at (0.5,0.4 ) {$0.78$};
          \node at (0.5,-0.4 ) {$0.52$};
          \node at (0.5,-1.2 ) {$0.26$};
          \node at (0.5,-2.0 ) {$0.0$};
    \end{tikzpicture}\\
    \centering
		\includegraphics[width=0.35\textwidth,trim={3cm 3cm 3cm 3cm},clip]{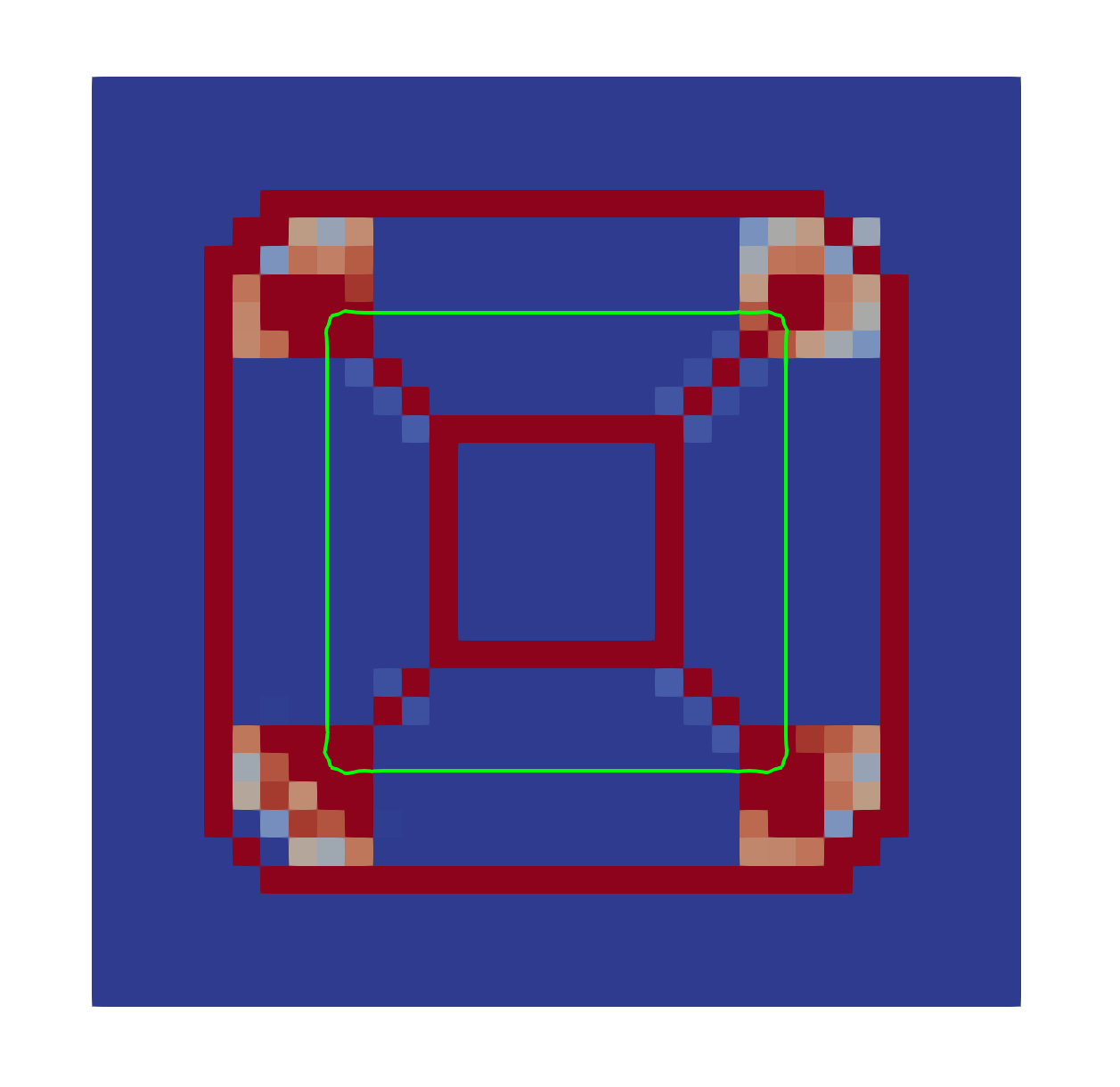}\hspace{0.5cm}
		\includegraphics[width=0.35\textwidth,trim={3cm 3cm 3cm 3cm},clip]{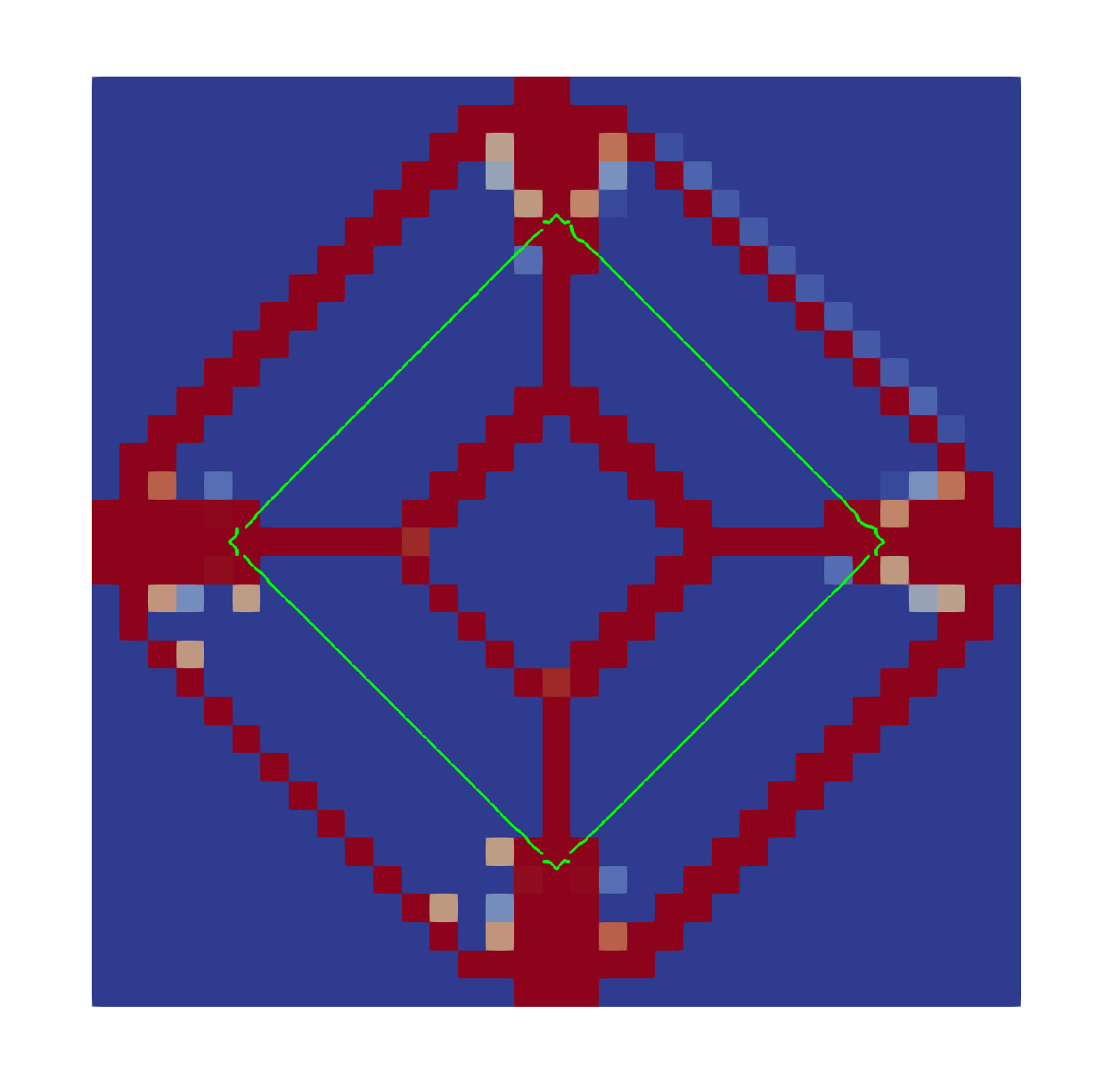}
    \begin{tikzpicture}[line cap=round,line join=round,>=triangle 45,x=1.15cm,y=1.15cm,axis/.style={->}]
          \node at (0 ,-0.4) {\includegraphics[scale=0.10]{./pictures/bar_color_rot.png}};
          \node at (0.0,1.6 ) {$\text{FV ratio}$};
          \node at (0.5,1.2 ) {$1.0$};
          \node at (0.5,0.4 ) {$0.75$};
          \node at (0.5,-0.4 ) {$0.5$};
          \node at (0.5,-1.2 ) {$0.25$};
          \node at (0.5,-2.0 ) {$0.0$};
    \end{tikzpicture}\\
    \caption{Absolute value of the level-set gradient (top row) and the amount of finite volume cells (bottom row) after $5000$ iterations. Initially a discontinuous rectangular shaped level-set is set. The level-set zero position is indicated by a green line.}
    \label{fig:rectangular}
\end{figure}
Figure \ref{fig:rectangular} shows the final results of the reinitialization procedure with the new regularization strategy at the steady state after $5000$ iterations with the third order Runge-Kutta method~\cite{Williamson1980}. The graph shows that the absolute value of the level-set gradient is very close to the desired one ($|\nabla \phi| \approx 1$), so the reinitialization is successful and a signed distance function is created. Moreover it is clearly visible that the FV sub-cell stabilization is only applied at the kinks of the level-set solution. This means that the solution is high order accurate in most of the vicinity of the level-set zero, since kinks at the level-set zero only occur at the corners of the squares.
\subsection{Reinitialization of strongly disturbed level-set function in 2d} \label{subsec:Hartmann}
%
This test case is taken from Russo and Smereka~\cite{russo2000} with the parameters according to Hartmann et al.~\cite{hartmann2008}. The test case illustrates the applicability of the introduced reinitialization strategy to a challenging problem with very strong gradients. The initial level-set function on the domain $\Omega=[-5,5]^2$ is given by
\begin{align}
  \phi(x,y)=g(x,y)\left( r- \sqrt{x^2 + y^2} \right),
\end{align}
where $g(x,y)$ is a disturbance function defined as
\begin{align}
  g(x,y)=0.1 +\left(x-r\right)^2+ \left(y-r\right)^2,
\end{align}
with the radius $r=3$ of the zero contour line. Two different meshes are used: a structured grid with $96^2$ elements and an unstructured grid with a comparable resolution of $9594$ elements. The polynomial degree is set to $\mathcal{N}=4$ and the remaining parameters for the simulation are: $\epsilon=20$, $CFL=0.9$, $\phi_\text{cut-off}=1$, upper threshold of indicator $\mathcal{S}_{up}=-5.5$ and lower threshold of indicator is $\mathcal{S}_{low}=-6.5$ with $n=2$.
\begin{figure}[h!t]
\begin{center}
	\begin{tikzpicture}[line cap=round,line join=round,>=triangle 45,x=1.15cm,y=1.15cm,axis/.style={->}]
        \node at (0 ,0) {\includegraphics[scale=0.16]{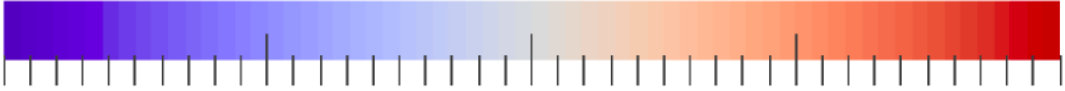}};
        \node at (3.5, 0.0 ) {$\text{FV ratio}$};
        \node at (2.6, -0.5 ) {$1.0$};
        \node at (1.3, -0.5 ) {$0.75$};
        \node at (0.0, -0.5 ) {$0.5$};
        \node at (-1.3,-0.5 ) {$0.25$};
        \node at (-2.6,-0.5 ) {$0.0$};
        \node at (-3.5, 0.0 ) {\color{white}{$\text{FV ratio}$}};
  \end{tikzpicture}\\
  $1$ iteration \hspace{3.3cm} $1$ iteration \vspace{0.2cm}\\
  \includegraphics[width=0.36\textwidth]{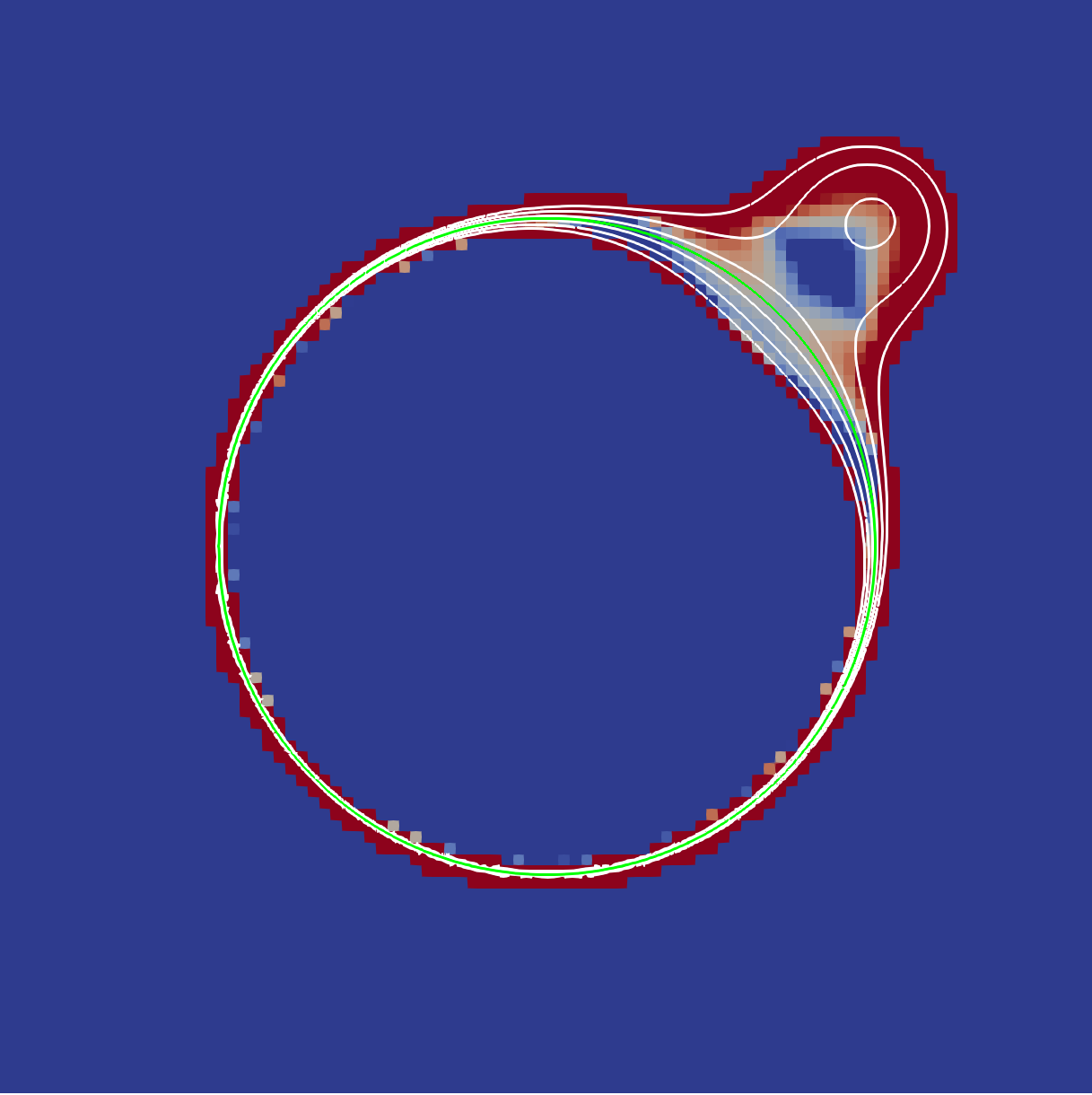} \hspace{0.5cm}
	\includegraphics[width=0.36\textwidth]{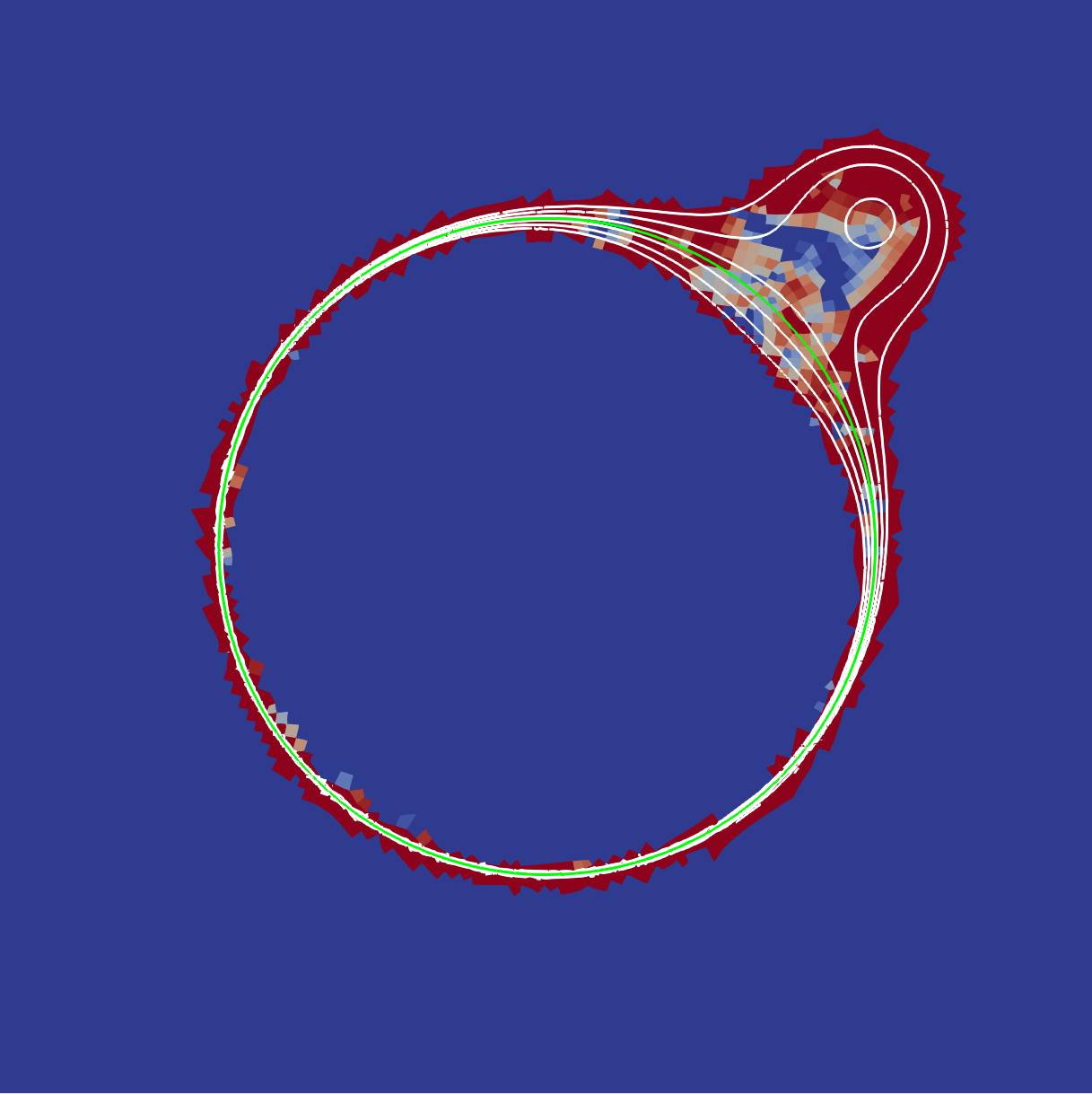} \\\vspace{-0.2cm}
	$200$ iterations  \hspace{2.7cm}  $600$ iterations \vspace{0.2cm}\\
  \includegraphics[width=0.36\textwidth]{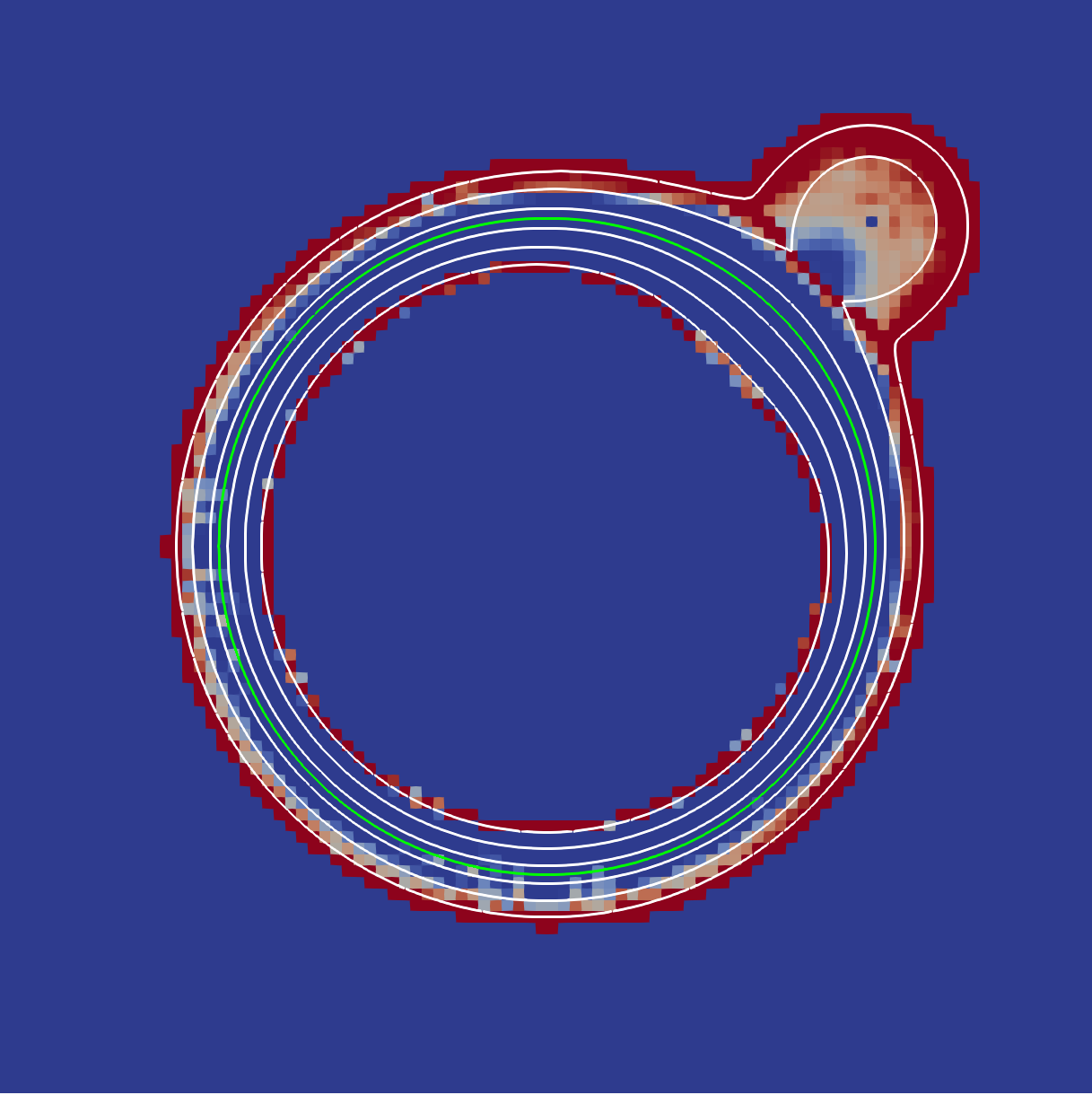}  \hspace{0.5cm}
	\includegraphics[width=0.36\textwidth]{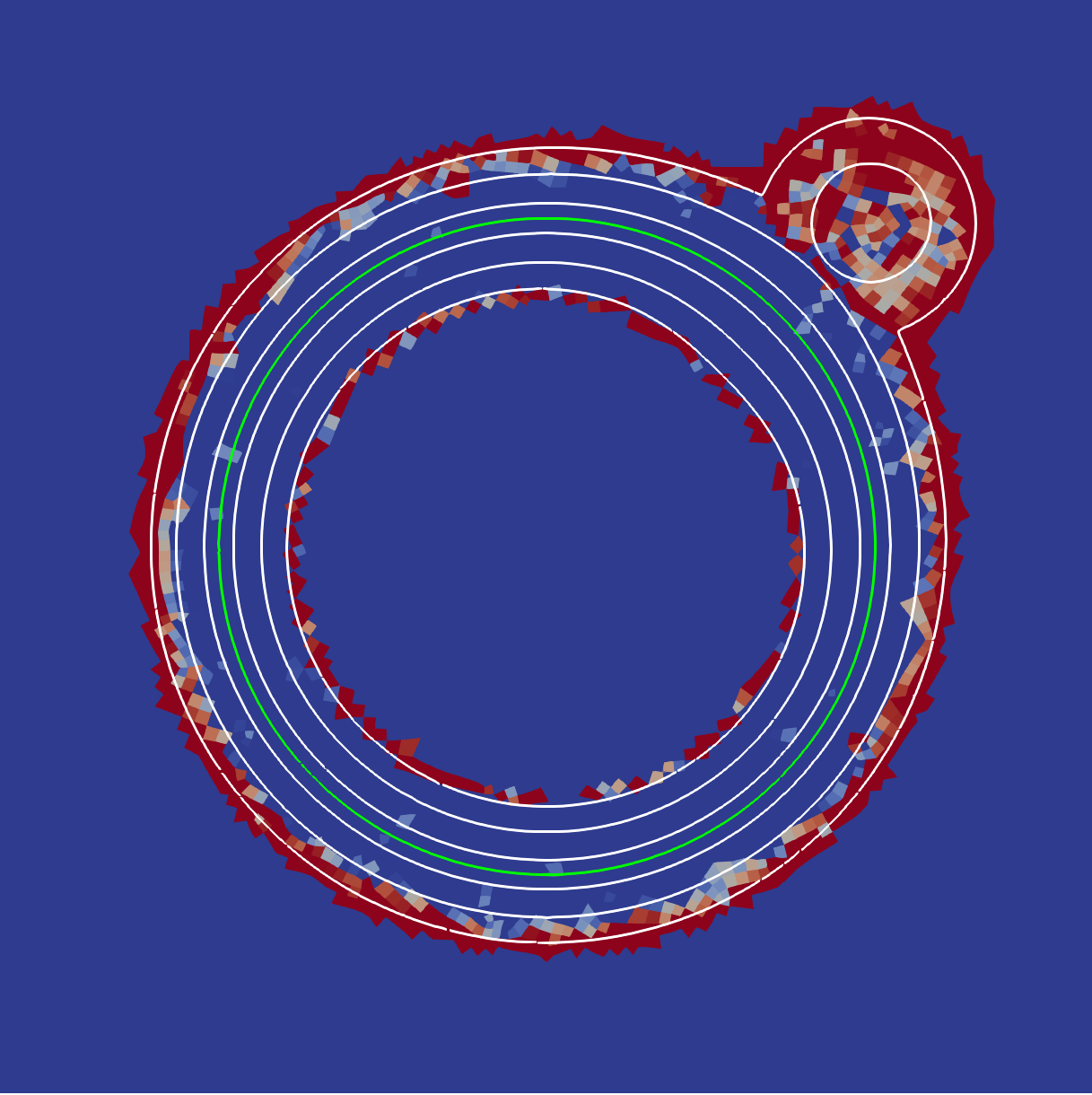} \\\vspace{-0.2cm}
	 $1000$ iterations  \hspace{2.5cm} $3000$ iterations \vspace{0.2cm}\\
  \includegraphics[width=0.36\textwidth]{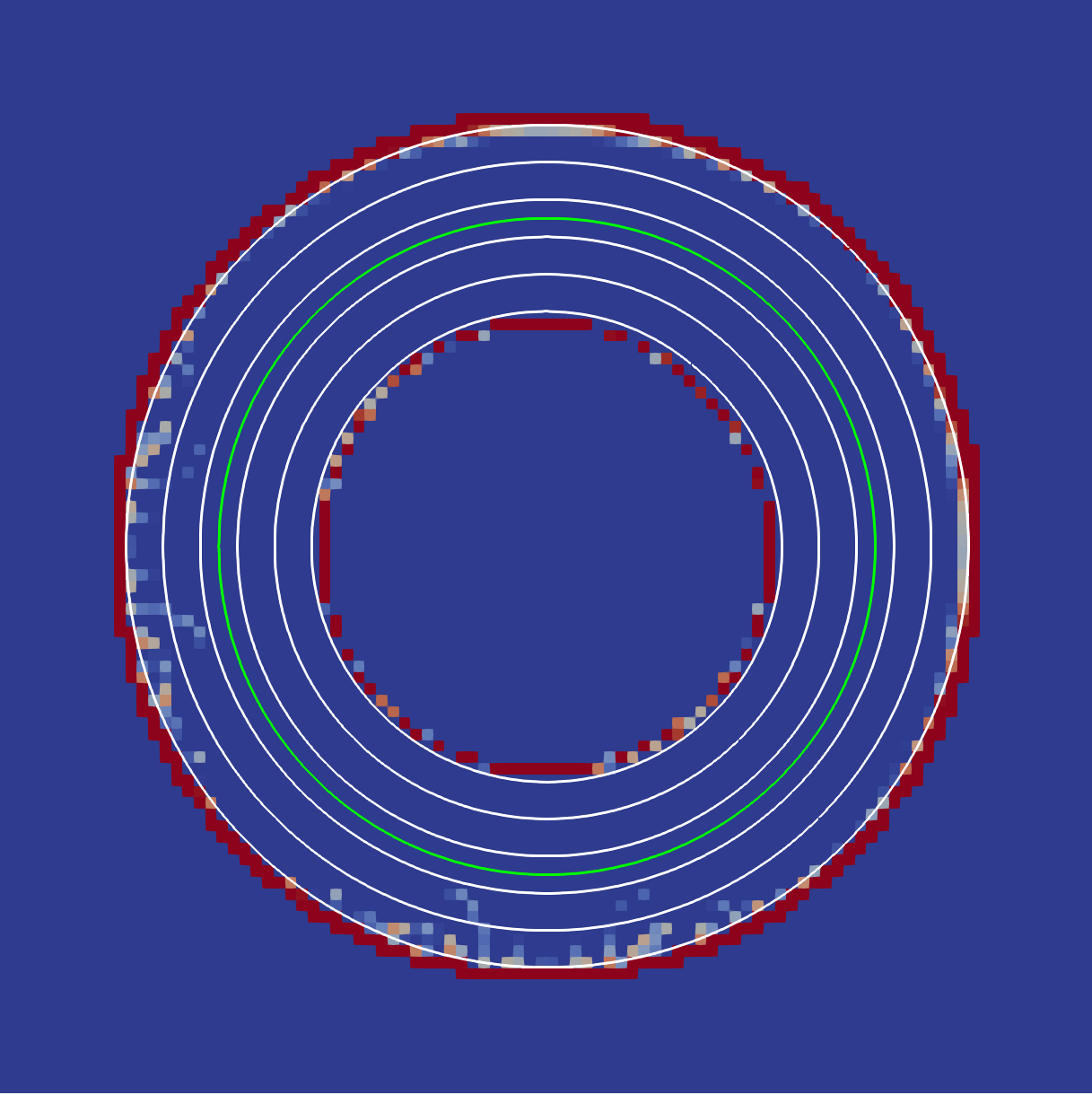} \hspace{0.5cm}
	\includegraphics[width=0.36\textwidth]{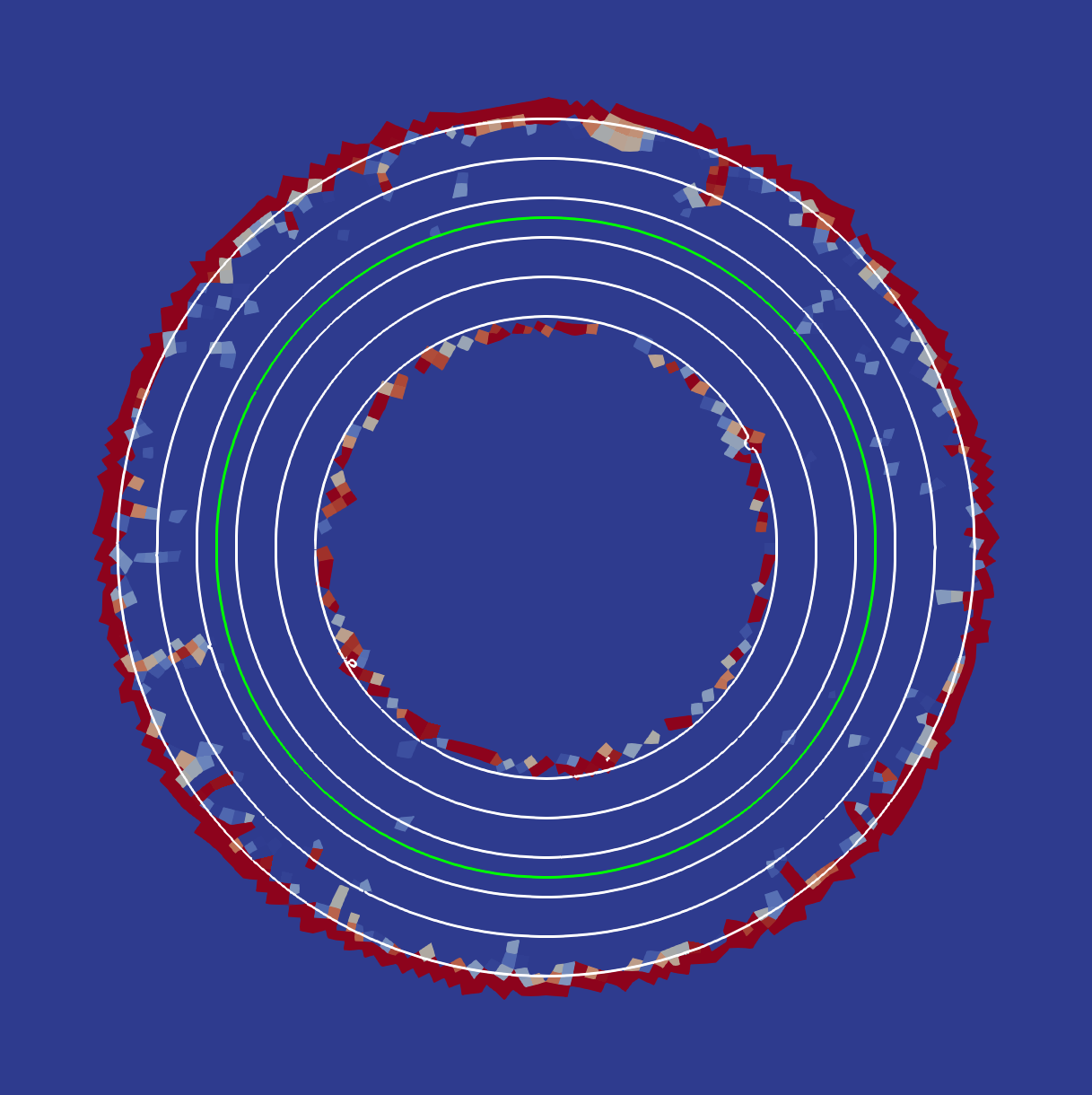} \\\vspace{-0.2cm}
	\caption{Evolution of level-set contours from $-0.9$ to $0.9$ and required amount of finite volume sub-cells for Hartmann test case for calculations on a Cartesian grid (left) and an unstructured grid (right). The level-set zero position is indicated by a green line.}
	\label{fig:hartmann}
\end{center}
\end{figure}
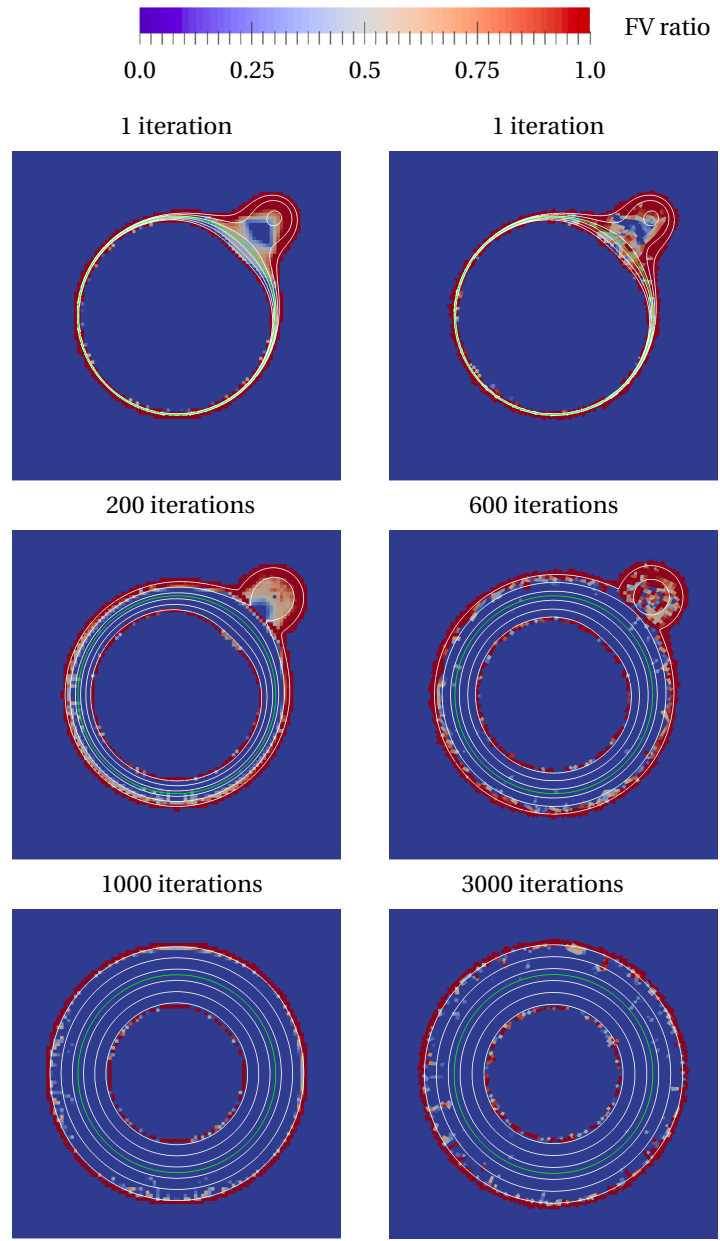
Figure~\ref{fig:hartmann} shows the isocontours of the level-set and the present amount of FV sub-cells in the computational domain for the structured and the unstructured mesh. Note that the time step of the third order Runge-Kutta method on the unstructured mesh is approximately three times smaller due to distorted elements. It is visible that the very strong gradient around the level-set zero position and the highly disturbed upper right corner require the use of stabilization. During the reinitialization procedure, the gradients are reduced and so less stabilization is necessary. Hence, the finite volume sub-cells in the vicinity of the level-set zero vanish. For the structured and the unstructured mesh a similar behavior is observed.
We can conclude that even for very difficult problems the converged solution requires no stabilization near the level-set zero positions and therefore provides a high order solution. This allows an accurate calculation of the curvature.
%
\subsection{Reinitialization of strongly disturbed level-set function in 3d} \label{subsec:3Dreinit}
\begin{figure}[h!t]
  \begin{tikzpicture}[line cap=round,line join=round,>=triangle 45,x=1.15cm,y=1.15cm,axis/.style={->}]
			\node at (0 ,0)   {\includegraphics[width=0.36\textwidth,trim={0 0 16cm 0},clip]{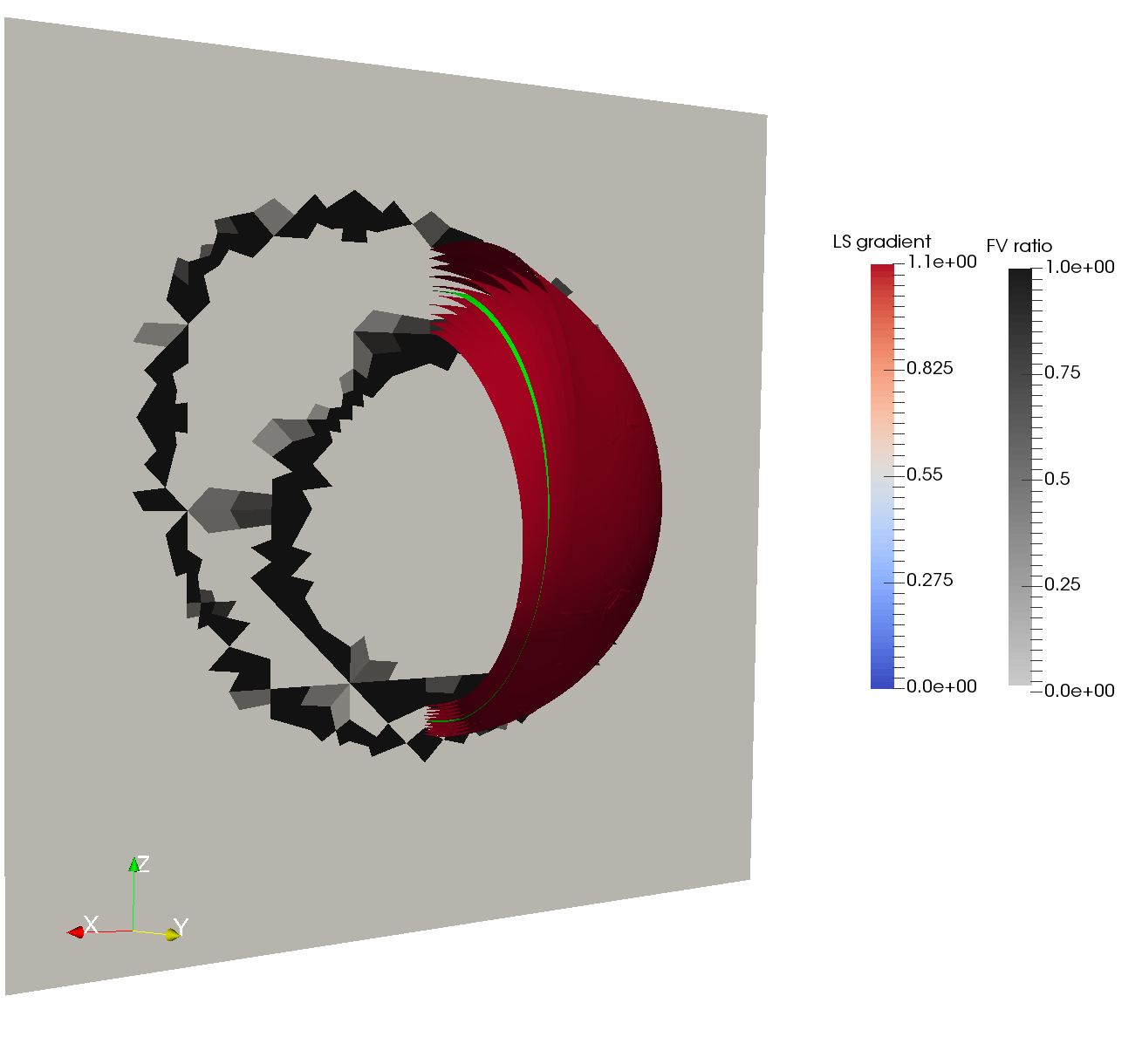}};
			\node at (3.4,-1) {\includegraphics[width=0.36\textwidth,trim={0 0 16cm 0},clip]{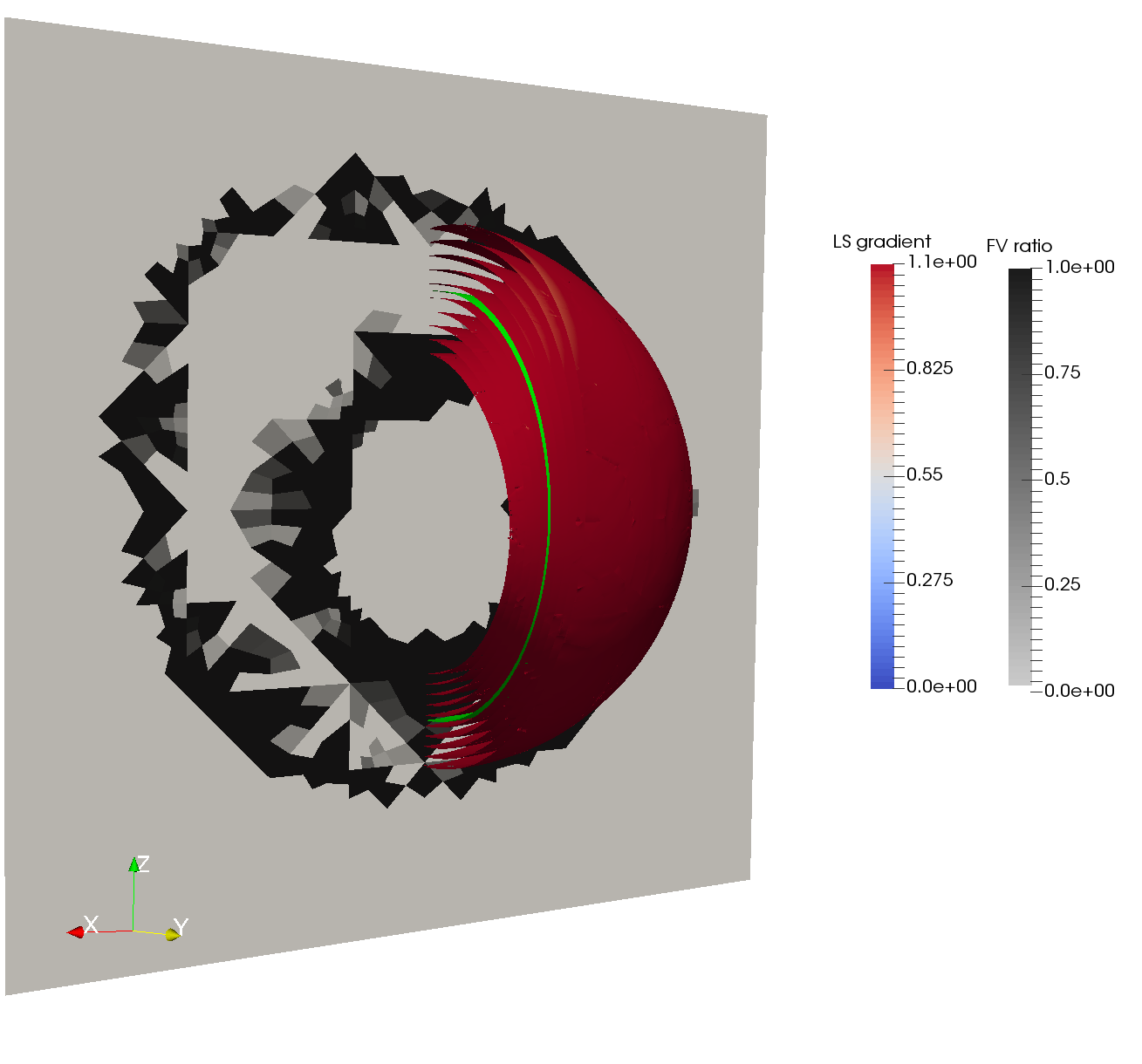}};
      \node at (7.0,-2) {\includegraphics[width=0.36\textwidth,trim={0 0 16cm 0},clip]{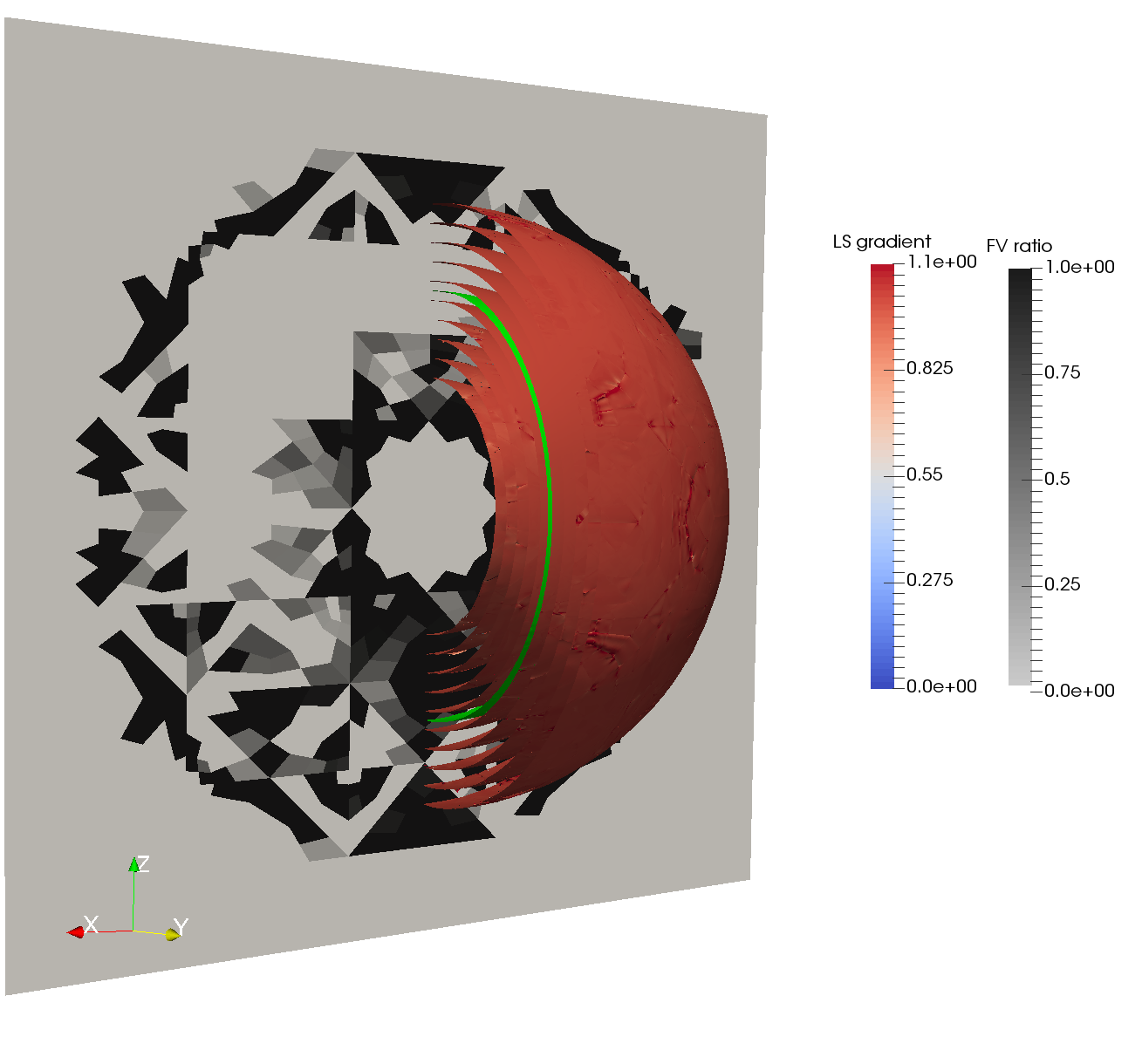}};
      \node at (6    , 2.3) {\includegraphics[scale=0.14]{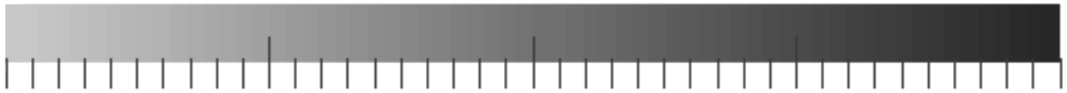}};
      \node at (6.0  , 2.8 ) {$\text{FV ratio}$};
      \node at (8.2  , 1.8 ) {$1.0$};
      \node at (7.1  , 1.8 ) {$0.75$};
      \node at (6.0  , 1.8 ) {$0.5$};
      \node at (4.9  , 1.8 ) {$0.25$};
      \node at (3.8  , 1.8 ) {$0.0$};
      \node at (1    , -4.5) {\includegraphics[scale=0.14]{./pictures/bar_color.png}};
      \node at (1.0  , -4.0 ) {$|\nabla \phi|$};
      \node at (3.2  , -5.0 ) {$1.1$};
      \node at (2.1  , -5.0 ) {$0.825$};
      \node at (1.0  , -5.0 ) {$0.55$};
      \node at (-0.1 , -5.0 ) {$0.275$};
      \node at (-1.2 , -5.0 ) {$0.0$};
      \node at (0.0 , -2.5 ) {$1$ iteration};
      \node at (3.5 , -3.5 ) {$300$ iterations};
      \node at (7.5 , -4.5 ) {$3000$ iterations};
  \end{tikzpicture}\\
  \caption{Evolution of level-set contours from $-0.4$ to $0.4$ colored by the level-set gradient. The required amount of finite volume sub-cells in the $y$-normal plane is indicated by the grey scale. The level-set zero position is indicated by the green isocontour.}\label{fig:hartmann_3d}
\end{figure}

To illustrate the applicability of the method to a three dimensional problem, we use the disturbed sphere described in~\cite{karakus2016755}. It is similar to the two-dimensional case in the previous section. On the domain $\Omega=[-2,2]^3$ the initial level-set function is given by
\begin{align}
	\phi(x,y,z)=((x-1)^2+(y-1)^2+(z-1)^2+0.1)\left(\sqrt{x^2 + y^2 + z^2}-1 \right).
\end{align}
A fully unstructured mesh with $96000$ elements and a polynomial degree of $\mathcal{N}=4$ are used. The mesh was generated from a Cartesian grid. The elements were split into tetrahedrons, which were split again to obtain hexahedrons. This procedure leads to strongly deformed grid cells, which cause numerical artifacts. They can be observed in the plot of the level-set gradient, Fig.~\ref{fig:hartmann_3d}.
The parameters for the simulation are: $\epsilon=20$, $CFL=0.9$ for $3^{rd}$ order Runge-Kutta method, $\phi_\text{cut-off}=0.6$, upper threshold of indicator $\mathcal{S}_{up}=-8.0$ and lower threshold of indicator $\mathcal{S}_{low}=-9.0$ with $n=1$ for Eq.~\eqref{eq:Indicator}. The gradients of the level-set are calculated with the BR1 lifting procedure and the central least squares method.
Figure~\ref{fig:hartmann_3d} shows, that the signed distance property of the level-set field is reached after $3000$ iterations, although the strongly disturbed level-set field has been initialized on a fully unstructured three dimensional mesh. Hence, we can conclude that the method allows the reinitialization on three dimensional unstructured meshes.

\section{Conclusion} \label{sec:Conclusion}
In this paper we developed a novel regularization strategy for the local discontinuous Galerkin method and applied it to reinitialize level-set functions. The novel approach is related to finite volume sub-cell shock-capturing concepts. An indicator is used to identify cells where the LDG scheme is unstable, e.g. at kinks and discontinuities and a first order FV operator on sub-cells is applied to capture them. This procedure stabilizes the numerical scheme. It is of arbitrary high order in the vicinity of the level-set zero, since the low order regularization is only applied at kinks/discontinuities, which typically occur away from the zero. Exceptions are regions where grid refinement is necessary or topological changes require the use of a low order discretization. The properties of the LDG scheme with respect to an easy parallelization and the use of unstructured grids are retained. Both points are benefits compared to finite volume or finite difference approaches. Compared with artificial viscosity approaches our method converges against the solution of the original equation. We also save the cost of evaluating a second order term and introduce the option to fall back to a minimal stencil, which might be useful if the method is applied to merging droplets or drop-wall interactions.
We presented h- and p-convergence studies of the introduced schemes and demonstrated the need for high order methods near the zero of the level-set to calculate the curvature. Afterwards we showed that our scheme can handle discontinuous initial conditions. We concluded the paper by solving challenging benchmark test problems to proof that our approach works well on both structured and unstructured grids in two and three space dimensions. In particular we demonstrated that even if a stabilization is required near the zero of the level-set it typically vanishes during the reinitialization process.
In the future, we plan to combine the presented algorithms with a ghost-fluid method~\cite{fechter2015} to evaluate its benefits compared with alternative approaches. The investigation of spurious currents as well as droplet collisions and drop-wall interactions are of particular interest.

\section*{Acknowledgments}

The authors kindly acknowledge the financial support provided by the German Research Foundation (DFG) through SFB-TRR 40 "Technological Foundations for the Design of Thermally and Mechanically Highly Loaded Components of Future Space Transportation Systems" and GRK 2160/1 "Droplet Interaction Technologies". 
The simulations were performed on the national supercomputer Cray XC40 (Hazel Hen) at the High Performance Computing Center Stuttgart (HLRS) under the grant number \emph{hpcmphas/44084} and the \emph{hpcdg} project.

\section*{References}

\bibliography{literature}

\end{document}